\documentclass[12pt]{amsart}
\usepackage{}

\usepackage{amsmath}
\usepackage{amsfonts}
\usepackage{amssymb}
\usepackage[all]{xy}           

\usepackage{bbding}
\usepackage{txfonts}
\usepackage{amscd}

\usepackage[shortlabels]{enumitem}
\usepackage{ifpdf}
\ifpdf
  \usepackage[colorlinks,final,backref=page,hyperindex]{hyperref}
\else
  \usepackage[colorlinks,final,backref=page,hyperindex]{hyperref}
\fi
\usepackage{tikz}
\usepackage[active]{srcltx}

\makeatletter

\newtheorem{thm}{Theorem}[section]
\newtheorem{lem}[thm]{Lemma}
\newtheorem{cor}[thm]{Corollary}
\newtheorem{pro}[thm]{Proposition}

\newtheorem{defi}[thm]{Definition}

\newcommand {\emptycomment}[1]{}

\newcommand{\cb}{\mathbb C}
\newcommand{\nb}{\mathbb N}

\newcommand{\Aut}{\mathrm{Aut}}

\newcommand{\CHom}{{\rm CHom}}

\def\Ker{\mathrm{Ker}\,}

\begin{document}

\title[The non-abelian extension and Wells map of Leibniz conformal algebra]
{The non-abelian extension and Wells map of Leibniz conformal algebra}

\author{Jun Zhao}
\address{School of Mathematics and Statistics, Henan University, Kaifeng 475004, China}
\email{zhaoj@henu.edu.cn}

\author{Bo Hou}
\address{School of Mathematics and Statistics, Henan University, Kaifeng 475004, China}
\email{bohou1981@163.com}

\author{Xin Zhou}
\address{School of Mathematics and Statistics, Yili Normal University, Yining 835000, China}
\email{zhoux566@nenu.edu.cn}
\vspace{-5mm}

\begin{abstract}
In this paper, we study the theory of non-abelian extensions of a Leibniz conformal algebra $R$ by a Leibniz conformal algebra $H$ and prove that all the non-abelian extensions are classified by non-abelian $2$nd cohomology $H^2_{nab}(R,H)$ in the sense of equivalence. Then we introduce a differential graded Lie algebra $\mathfrak{L}$ and show that the set of its Maurer-Cartan elements in bijection with the set of non-abelian extensions. Finally, as an application of non-abelian extension, we consider the inducibility of a pair of automorphisms about a non-abelian extension, and give the fundamental sequence of Wells of Leibniz conformal algebra $R$. Especially, we discuss the extensibility problem of derivations about an abelian extension of $R$.
\end{abstract}


\maketitle




\vspace{-4mm}

\section{Introduction}\label{sec:intr}

The notion of Lie conformal algebra encodes an axiomatic description of the operator product expansion of chiral fields in two-dimensional conformal field theory. It has been proved to be an effective tool for the study of infinite-dimensional Lie algebras satisfying the locality property \cite{Ka,Ka1}. And it is also closely related to vertex algebra \cite{Ka}. As a generalization of Lie conformal algebra, the study of Leibniz conformal algebra has received a lot of attention.
The cohomology theory of Leibniz conformal algebra were studied in \cite{Z}. The references \cite{DS} and \cite{GW} studied homotopification and categorification of Leibniz conformal algebra and the twisted relative Rota-Baxter operators on Leibniz conformal algebra respectively. The authors in \cite{HY} studied the extension structures of Leibniz conformal algebra. The Leibniz conformal algebras of rank two and rank three were studied by Z. Wu in \cite{Wu} and \cite{Wu1} respectively. The authors in \cite{XBY} discussed about Leibniz conformal bialgebras and the classical Leibniz conformal Yang-Baxter equation of Leibniz conformal algebra. And the authors in  also \cite{ZH} considered the quadratic Leibniz conformal algebras.

In recent years, extension problem, especially non-abelian extension, plays a very important role in algebra theory and it is closely related to representation and cohomology theory. And it has attracted much attention due to its connections with group \cite{MOS,W}, Lie algebra \cite{F}, pre-Lie algebra \cite{SZ,MSW}, associative algebra \cite{A}, $3$-Lie algebra  \cite{TX}, Rota-Baxter Lie algebra \cite{MDH}, Rota-Baxter Leibniz algebra \cite{GH}, $3$-Leibniz algebra \cite{XS}, associative conformal algebra \cite{HZ}, Nijenhuis Lie algebras \cite{D}, Lie H-pseudoalgebras \cite{D1} and so on. In addition, the inducibility of a pair of automorphisms about a non-abelian extension and the extensibility problem of derivations about an abelian extension are also interest directions that scholars of mathematics are interested in. In this paper, we study the non-abelian extensions of a Leibniz conformal algebra. And as applications of extension problem, we discuss the inducibility of a pair of automorphisms and the extensibility problem of derivations.

The paper is organized as follows. In Section 2, we recall some related definitions, such as differential graded Lie algebra, Leibniz conformal algebra, Maurer-Cartan element and the cohomology of a Leibniz conformal algebra $R$ and so on. In Section 3, firstly, we give the definitions of the non-abelian extension and non-abelian $2$-cocycle of a Leibniz conformal algebra $R$. And then we prove that all non-abelian extensions of Leibniz conformal algebra $R$ by a Leibniz conformal algebra $H$ are classified by $2$nd cohomology space $H^2_{nab}(R,H)$. In Section 4, given two Leibniz conformal algebras $R$ and $H$, we construct a differential graded Lie algebra $\mathcal{L}$ and show that the set of its Maurer-Cartan elements in bijection with the set of extensions of $R$ by $H$, which gives a interpretation of $H^2_{nab}(R,H)$. In Section 5, we study the inducibility of a pair of automorphisms about a non-abelian extension of Leibniz conformal algebra $R$ by a Leibniz conformal algebra $H$. In Section 6, we consider the extensibility problem of derivations about an abelian extension of Leibniz conformal algebra $R$ by a $\cb[\partial]$-module $H$.

Throughout this paper, all vector spaces, linear maps, and tensor products are over the complex field $\mathbb{C}$. Use $\mathbb{Z}$ to denote the set of integers.

\bigskip
\section{The preliminaries of a Leibniz conformal algebra}
In this section, we recall some relative definitions about differential graded Lie algebra and Leibniz conformal algebra. The following notions can be found in \cite{F} and \cite{ZH}.

\begin{defi}
A graded Lie algebra is a $\mathbb{Z}$-graded vector space $L=\oplus_{n\in\mathbb{Z}}L_n$ equipped with a degree preserving bilinear bracket $[\cdot,\cdot]:L\times L\rightarrow L$,
which satisfies
\begin{enumerate}
\item[$(G1)$] graded antisymmetry: $ [a,b]=-(-1)^{|a||b|}[b, a],$
\item[$(G2)$] graded Leibniz rule:  $[a,[b,c]]=[[a,b],c]+(-1)^{|a||b|}[b,[a,c]].$
\end{enumerate}
Here $a,b,c$ are homogeneous elements of $L$ and the degree $|x|$ of a homogeneous element $x\in L_n$ is defined by $n$.
\end{defi}

\begin{defi}
A differential graded Lie algebra is a graded Lie algebra $(L,[\cdot,\cdot])$ equipped with a homological derivation $d:L\rightarrow L$ of degree $1$,
which satisfies
\begin{enumerate}
\item[$(D1)$]  $|da|=|a|+1$ ($d$ of\ degree\ $1$),
\item[$(D2)$]  $d([a,b])=[d(a),b]+(-1)^{|a|}[a,d(b)]$ (derivation),
\item[$(D3)$]  $d^2=0$ (homological).
\end{enumerate}
\end{defi}

\begin{defi}
The set $MC(L)$ of Maurer-Cartan elements of a differential graded Lie algebra $(L,[\cdot,\cdot],d)$ is defined to be
$$\{\alpha\in L^1\mid d\alpha+\frac{1}{2}[\alpha,\alpha]=0\}.$$
\end{defi}

We now recall the equivalence relation on $MC(L)$. Let $\beta$ and $\beta'$ be two elements in $MC(L)$. We say that $\beta$ and $\bar{\beta}$ are equivalent if there exists an $\alpha\in L^0$ such that
$$\beta'=e^{\rm ad\alpha}(\beta)-\frac{e^{\rm ad\alpha}-1}{\rm ad\alpha}d(\alpha).$$
The set of the equivalence relation on $MC(L)$ is denoted by $\mathcal{MC}(L)$.

\begin{defi}
A Leibniz conformal algebra is a $\mathbb{C}[\partial]$-module $R$ endowed
with a bilinear map $R\times R\rightarrow R[\lambda], (a,b)\mapsto[a_\lambda b]$ for $a,b\in R$, called the $\lambda$-bracket, satisfying the
following axioms ($c\in R$):
\begin{enumerate}
\item[$(L1)$] Conformal sesquilinearity: $ [(\partial a)_\lambda b]=-\lambda[a_\lambda b], \ \ [a_\lambda \partial b]=(\partial +\lambda)[a_\lambda b];$
\item[$(L2)$] Jacobi identity: $[a_\lambda[b_\mu c]]=[[a_\lambda b]_{\lambda+\mu}c]+[b_\mu[a_\lambda c]].$
\end{enumerate}
\end{defi}
A Leibniz conformal algebra $R$ is called {\it finite} if it is finitely generated as a $\cb[\partial]$-module. The {\it rank} of a Leibniz conformal algebra $R$ is its rank as a $\cb[\partial]$-module. This paper we assume that all $\cb[\partial]$-modules are finitely generated.

\begin{defi}
Let $M_1$ and $M_2$ be two $\cb[\partial]$-modules. A {\rm conformal linear map} from $M_1$ to $M_2$ is a $\cb$-linear map $f:M_1\rightarrow M_2[\lambda],$ denoted by $f_\lambda:M_1\rightarrow M_2[\lambda]$, such that $f_\lambda(\partial v)=(\partial+\lambda)f_\lambda(v)$ for $v\in M_1$, usually written as $f$.

Moreover, let $N$ also be a $\cb[\partial]$-module. A {\rm conformal bilinear map} from $M_1\times M_2$ to $N$ is a
$\cb$-bilinear map $f:M_1\times M_2\rightarrow N[\lambda]$, denoted by $f_\lambda:M_1\times M_2\rightarrow N[\lambda]$, such that $f_\lambda(\partial u,v)=-\lambda f_\lambda(u,v)$ and $f_\lambda(u,\partial v)=(\partial+\lambda)f_\lambda(u,v)$ for any $u\in M_1, v\in M_2$.
\end{defi}

\begin{defi}\label{derivation}
A left derivation on Leibniz conformal algebra $R$ is a $\cb[\partial]$-module homomorphism $d:R\rightarrow R$ satisfying
\begin{eqnarray*}d([a_\lambda b])=[a_\lambda d(b)]+[{d(a)}_{\lambda} b],\qquad \forall \ a,b\in R.
\end{eqnarray*}
\end{defi}

\begin{defi}
A right derivation on Leibniz conformal algebra $R$ is $\cb[\partial]$-module homomorphism $d:R\rightarrow R[\lambda]$ satisfying
\begin{eqnarray*}d([a_\lambda b])=[a_\lambda d(b)]-[b_{-\partial-\lambda}d(a)],\qquad \forall \ a,b\in R.
\end{eqnarray*}
\end{defi}
Denote by ${\rm Der}^L(R)$ and ${\rm Der}^R(R)$ the set of all left derivations and the set of all right derivations on Leibniz conformal algebra $R$ respectively. Then we can obtain that
${\rm Der}^L(R)$ and ${\rm Der}^R(R)$ are Lie algebras.

\begin{defi}
A left conformal derivation on Leibniz conformal algebra $R$ is a a conformal linear map $d:R\rightarrow R$ satisfying
\begin{eqnarray*}d_\lambda([a_\mu b])=[a_\mu d_\lambda(b)]+[{d_\lambda(a)}_{\lambda+\mu} b],\qquad \forall \ a,b\in R.
\end{eqnarray*}
\end{defi}

\begin{defi}
A right conformal derivation on Leibniz conformal algebra $R$ is a conformal linear map $d:R\rightarrow R[\lambda]$ satisfying
\begin{eqnarray*}d_\lambda([a_\mu b])=[a_\mu d_\lambda(b)]-[b_{-\partial-\mu}d_\lambda(a)],\qquad \forall \ a,b\in R.
\end{eqnarray*}
\end{defi}
Denote by ${\rm CDer}^L(R)$ and ${\rm CDer}^R(R)$ the set of all left conformal derivations and the set of all right conformal derivations on Leibniz conformal algebra $R$ respectively. Then we can obtain that
${\rm CDer}^L(R)$ and ${\rm CDer}^R(R)$ are Lie conformal algebras. Furthermore, if we denote $({\rm ad}^L_a)_\lambda(b)=[a_\lambda b]$ and $({\rm ad}^R_a)_\lambda (b)=[b_{-\partial-\lambda} a]$, then ${\rm ad}^L_a$ is a left conformal derivation and ${\rm ad}^R_a$ is a right conformal derivation on $R$.

\begin{defi}
A homomorphism of Leibniz conformal algebras from $R$ to $R'$ is a $\cb[\partial]$-module homomorphism $\phi:R\rightarrow R'$ such that the following condition is satisfied:
\begin{eqnarray}\phi(a_\lambda b)=\phi(a)_\lambda\phi(b),\qquad \forall \ a,b\in R.\label{defi2.6}
\end{eqnarray}
\end{defi}

\begin{defi}
A conformal module $M$ over a Leibniz conformal algebra $R$ is a $\cb[\partial]$-module with two $\cb$-linear maps $R\otimes M\rightarrow M[\lambda], a\otimes v\mapsto a_\lambda v$, $M\otimes R\rightarrow M[\lambda], v\otimes a\mapsto v_\lambda a$, such that
\begin{eqnarray*}
&&(\partial a)_\lambda v=-\lambda a_\lambda v,\qquad a_\lambda (\partial v)=(\partial+\lambda)a_\lambda v,\\
&&(\partial v)_\lambda a=-\lambda v_\lambda a,\qquad v_\lambda (\partial a)=(\partial+\lambda)v_\lambda a,\\
&&a_\lambda(b_\mu v)=[a_\lambda b]_{\lambda+\mu}v+b_\mu(a_\lambda v),\\
&&a_\lambda(v_\mu b)=(a_\lambda v)_{\lambda+\mu}b+v_\mu[a_\lambda b],\\
&&v_\lambda[a_\mu b]=(v_\lambda a)_{\lambda+\mu}b+a_\mu(v_\lambda b).
\end{eqnarray*}
hold for $a,b\in R$ and $v\in M$.
\end{defi}

Let $R$ be a Leibniz conformal algebra. Define the module action of $R$ on itself is $a_\lambda b:=[a_\lambda b]$ for $a, b\in R$. Then $R$ is a module on itself. We call this module {\it adjoint module}.

\begin{pro}
Let $R$ be a Leibniz conformal algebra and $M$ be an $R$-module. Then the $\cb[\partial]$-module $R\oplus M$ is a Leibniz conformal algebra with the following $\lambda$-product:
$$[(a+u)_\lambda(b+v)]=[a_\lambda b]+a_\lambda v+u_\lambda b,\qquad \forall\ a,b\in R, u,v\in M.$$
We call it {\rm semi-direct product Leibniz conformal algebra} of $R$ and $M$, and denote it by $R\ltimes M$.
\end{pro}

Let $R$ and $M$ be two $\cb[\partial]$-modules. Denote by $\CHom(R^{\otimes n},M)$ the set of $\cb$-linear maps $\gamma:R^{\otimes n}\rightarrow M$ satisfying conformal antilinearity: \begin{eqnarray*}&&\gamma_{\lambda_{1},\cdots,\lambda_{n-1}}(a_1,\cdots,\partial a_i,\cdots,a_n) =-\lambda_{i}\gamma_{\lambda_{1},\cdots,\lambda_{n-1}}(a_1,\cdots,a_i,\cdots,a_n), \  i=1,\ldots,n-1,\label{ant1}\\
&&\gamma_{\lambda_{1},\cdots,\lambda_{n-1}}(a_1,\cdots, \partial a_n) =(\partial+\lambda_1+\ldots+\lambda_{n-1})\gamma_{\lambda_{1},\cdots,\lambda_{n-1}}(a_1,\cdots,a_n).\label{ant2}\end{eqnarray*}

Next, we recall the cohomology of a Leibniz conformal algebra $R$ in \cite{Z}. An $n$-cochain $(n\geq1)$ of the Leibniz conformal algebra $R$ with coefficients in an $R$-module $M$ is a $\cb$-linear map in $\CHom(R^{\otimes n},M)$. Denote by $C^n(R,M):=\CHom(R^{\otimes n},M)$ the set of all $n$-cochains of $R$.  For $\gamma\in C^n(R,M)$, $a_1,\ldots,a_{n+1}\in R$, the coboundary operator $\delta:C^n(R,M)\rightarrow C^{n+1}(R,M)$ is given by
\begin{eqnarray*}
&&(\delta \gamma)_{\lambda_1,\ldots,\lambda_n}(a_1,\ldots,a_{n+1})\\
&=&\sum\limits^n_{i=1}(-1)^{i+1}{a_i}_{\lambda_i}\gamma_{\lambda_1,\ldots,\hat{\lambda}_i,\ldots,\lambda_n}(a_1,\ldots,\hat{a}_i,\ldots,a_{n+1})+(-1)^{n+1} \gamma_{\lambda_1,\ldots,\lambda_{n-1}}
(a_1,\ldots,a_{n})_{\lambda_1+\cdots+\lambda_{n}}a_{n+1}\\
&&+\sum\limits_{1\leq i<j\leq n+1}(-1)^{i}\gamma_{\lambda_1,\ldots,\hat{\lambda}_{i},\ldots,\lambda_{j-1},\lambda_i+\lambda_j,\ldots,\lambda_{n-1}}
(a_1,\ldots,\hat{a}_{i},\ldots,a_{j-1},[{a_i}_{\lambda_i}a_j],\ldots,a_{n+1}).
\end{eqnarray*}

\begin{lem}
The map $\delta$ satisfies $\delta^2=0.$
\end{lem}
The cohomology space of the complex $(C^{\bullet}(R,M),\delta)$ will be denoted by $H^{\bullet}(R,M)$. It is obvious that $C^1(R,M)$ consists of all $\cb[\partial]$-module homomorphisms from $R$ to $M$ and $Z^1(R,M)$ consists of left derivations of from $R$ to $M$. In particular, if $R$ is a Leibniz conformal algebra and view $R$ as a conformal module on $R$, we denote the corresponding complex and cohomology space by $C^{\bullet}(R,R)$ and $H^{\bullet}(R,R)$. Moreover, $C^{\bullet+1}(R,R)=\oplus_{n\geq0}C^{n+1}(R,R)$ comes equipped with the Nijenhuis-Richardson graded Lie bracket
$$[\gamma,\theta]=\gamma\cdot\theta-(-1)^{mn}\theta\cdot\gamma,$$
for any $\gamma\in C^{m+1}(R,R), \theta\in C^{n+1}(R,R)$, where
\begin{eqnarray*}
&&(\gamma\cdot\theta)_{\lambda_1,\cdots,\lambda_{m+n}}(a_1,\cdots,a_{m+n+1})\\
&=&\sum_{0\leq k\leq m}(-1)^{k(n+1)}\sum\limits_\sigma {\rm sign}(\sigma)\gamma_{\lambda_{\sigma(1)},\cdots,\lambda_{\sigma(k)},\lambda_{\sigma(k+1)}+\cdots+\lambda_{\sigma(k+m)}+\lambda_{k+m+1},\lambda_{k+m+2},\cdots,\lambda_{n+m}}\\
&&\big(a_{\sigma(1)},\cdots,a_{\sigma(k)},\theta_{\lambda_{\sigma(k+1)},\cdots,\lambda_{\sigma(k+m)}}(a_{\sigma(k+1)},,\cdots,a_{\sigma(k+m)},a_{k+m+1}),\cdots,a_{n+m+1}\big),
\end{eqnarray*}
with $\sigma\in S_{(k,n)}$, the set of permutations such that $\sigma(1)<\sigma(2)<\cdots<\sigma(k)$ and $\sigma(k+1)<\sigma(k+2)<\cdots<\sigma(k+n)$. Thus, $\big(C^{\bullet+1}(R,R),[\cdot,\cdot]\big)$ is a graded Lie algebra. Furthermore, the triple $\big(C^{\bullet+1}(R,R),[\cdot,\cdot],\delta\big)$ is a differential graded Lie algebra.

\bigskip

\section{Extensions of Leibniz conformal algebras classified by non-abelian cohomology}

In this section, we define the non-abelian extension and non-abelian $2$-cocycle of a Leibniz conformal algebra $R$ and obtain that all non-abelian extensions of Leibniz conformal algebra $R$ by a Leibniz conformal algebra $H$ are classified by $2$nd cohomology space $H^2_{nab}(R,H)$.
\begin{defi}
Let $R$ and $H$ be two Leibniz conformal algebras. An {\rm extension} $E$ of $R$ by $H$ is a short exact sequence of the form
$$0 \to H \stackrel{\alpha}\to E \stackrel{p}\to R \to 0.$$
\end{defi}

\begin{defi}
Let $E$ and $E'$ be two extensions of $R$ by $H$. They are said to be {\rm equivalent} if there exists a commutative diagram
$$\xymatrix
{0\ar[r]& H\ar[r]\ar@{=}[d]& E\ar[r]\ar[d]^{\varphi}& R\ar[r]\ar@{=}[d]& 0\\
0\ar[r]& H\ar[r]& E'\ar[r]& R\ar[r]& 0,}$$
where $\varphi:E\rightarrow E'$ is a Leibniz conformal algebra homomorphism.
\end{defi}

\begin{defi}\label{2-cocycle}
A non-abelian $2$-cocycle on $R$ with values in $H$ is a triple $(l,r,\chi)$, where $\chi:R\times R\rightarrow H[\lambda]$ is a conformal bilinear map and $l:R\rightarrow {\rm CDer}^L(H)$, $r:R\rightarrow {\rm CDer}^R(H)$ are $\cb[\partial]$-module homomorphisms satisfying
\begin{eqnarray*}
&&[(l(a)_\lambda h_1+r(a)_\lambda h_1)_{\lambda+\mu}h_2]=0,\ \ \ \ r(a)_\lambda(r(b)_\mu h+l(b)_\mu h)=0,\\
&&[{l(a)}_\lambda l(b)]-l([a_\lambda b])={\rm ad}^L_{\chi_\lambda(a,b)}\label{abel1},\\
&&[{l(a)}_\lambda r(b)]-r([a_\lambda b])={\rm ad}^R_{\chi_\lambda(a,b)}\label{abel1},\\
&&\chi_\lambda(a,[b_\mu c])-\chi_\mu(b,[a_\lambda c])-\chi_{\lambda+\mu}([a_\lambda b],c)\\
&&+{l(a)}_\lambda\chi_\mu(b,c)-{l(b)}_\mu\chi_\lambda(a,c)+{r(c)}_{-\partial-\lambda-\mu}\chi_\lambda(a,b)=0.\nonumber
\end{eqnarray*}
for any $a,b,c\in R, h_1,h_2\in H$.
\end{defi}
A {\it section} $s$ of the extension $E$ is a $\cb[\partial]$-module homomorphism $s:R\rightarrow E$ satisfying $p\circ s={\rm id}_R$.
\begin{defi}
Two cocycles $(l,r,\chi)$ and $(l',r',\chi')$ are said to be {\rm equivalent} if there exists a $\cb[\partial]$-module homomorphism $\beta:R\rightarrow H$ satisfying
\begin{eqnarray}
&&l'(a)-l(a)={\rm ad}^L_{\beta(a)},\ \ \ \ r'(a)-l(a)={\rm ad}^R_{\beta(a)},\label{cocycle1}\\
&&\chi'_\lambda(a,b)=\chi_\lambda(a,b)+{l(a)}_\lambda\beta(b)+{r(b)}_{-\partial-\lambda}\beta(a)-\beta([a_\lambda b])+[{\beta(a)}_\lambda\beta(b)].\label{cocycle2}
\end{eqnarray}
Denote it by $(l,r,\chi)\sim (l',r',\chi')$.
\end{defi}

We denote by $Z^2_{nab}(R,H)$ the set of non-abelian cocycles. Non-abelian cohomology $H^2_{nab}(R,H)$ will be the quotient of $Z^2_{nab}(R,H)$ by the equivalence relation.

\begin{pro}\label{class}
Let $R$ and $H$ be two Leibniz conformal algebras. Then all extensions of $R$ by $H$ are classified by $H^2_{nab}(R,H)$.
\end{pro}
\begin{proof}
Given an extension $E$ of $R$ by $H$, there exists a section $s:R\rightarrow E$ in
$$0 \to H \stackrel{\alpha}\to E \stackrel{p}\to R \to 0$$
such that $p\circ s={\rm id}_R$. The associated cocycle $(l,r,\chi)$ is defined by the formulas
\begin{eqnarray}
{l(a)}_\lambda(h):=[{s(a)}_\lambda h], \ \ \ {r(a)}_\lambda(h):=[h_{-\partial-\lambda}s(a)],\ \ \ {\chi_\lambda}(a,b):=[{s(a)}_\lambda s(b)]-s([a_\lambda b]),\label{extension}
\end{eqnarray}
for any $a,b\in R, h\in H$. Then $(l,r,\chi)\in Z^2_{nab}(R,H)$ is independent of the choice of the section $s$. The details as follows, let $s'$ be another section. For any $a\in R$, we have $p(s'(a)-s(a))=0$. Then $s'(a)-s(a)\in {\rm Ker} p={\rm Im} \alpha=H$. Define $\beta:R\rightarrow H$ as $\beta(a)=s(a)-s'(a)$. Hence
\begin{eqnarray*}
&&{(l(a)-l'(a))}_\lambda(h)=[({s(a)-s'(a))}_\lambda h]={({\rm ad}^L_{\beta(a)})}_\lambda h,\\
&&{(r(a)-r'(a))}_\lambda(h)=[h_{-\partial-\lambda}({s(a)-s'(a))}]={({\rm ad}^R_{\beta(a)})}_\lambda h,\\
&&{(\chi_\lambda-{\chi'}_\lambda)}(a,b)\\
&=&[{s(a)}_\lambda s(b)]-s([a_\lambda b])-[{s'(a)}_\lambda s'(b)]+s'([a_\lambda b])\\
&=&[{(s(a)-s'(a))}_\lambda (s(b)-s'(b))]+[{(s(a)-s'(a))}_\lambda s'(b)]\\
&&+[{s'(a)}_\lambda(s(b)-s'(b))]-(s-s')([a_\lambda b])\\
&=&{l'(a)}_\lambda \beta(b)+{r'(b)}_{-\partial-\lambda} \beta(a)-\beta([a_\lambda b])+[{\beta(a)}_\lambda\beta(b)].
\end{eqnarray*}
That is $(l,r,\psi)$ and $(l',r',\psi')$ are equivalent.

Suppose that two extensions $E$ and $E'$ of $R$ by $H$ are equivalent, i.e. there exists a commutative diagram
$$\xymatrix
{0\ar[r]& H\ar[r]\ar@{=}[d]& E\ar[r]\ar[d]^{\varphi}& R\ar[r]\ar@{=}[d]& 0\\
0\ar[r]& H\ar[r]& E'\ar[r]& R\ar[r]& 0,}$$
where $\varphi:E\rightarrow E'$ is a Leibniz conformal algebra homomorphism.
Choose two sections $s:R\rightarrow E$ and $\tilde{s}:R\rightarrow E'$. By direct computation, we have
$$(l,r,\chi)\sim (\tilde{l},\tilde{r},\tilde{\chi})$$
by the linear map $\beta_\varphi:=\varphi^{-1}\circ \tilde{s}-s$.

Conversely, for any a cocycle $(\psi,\chi)$, we can define an extension by
$$0 \to H \stackrel{}\to R\oplus H \stackrel{}\to R \to 0.$$
The $\lambda$-bracket defined on $R\oplus H$ is
$$[(a,u)_\lambda(b,v)]=\big([a_\lambda b],l(a)_\lambda(v)+r(b)_{-\partial-\lambda}(u)+\chi_\lambda(a,b)+[u_\lambda v]\big).$$
Note that the cocycle induced by the extension is in consistent with the given one $(l,r,\chi)$, where the section can be given by $s:R\rightarrow R\oplus H, a\mapsto (a,0)$.

Furthermore, if two cocycles $(l,r,\chi)$ and cocycles $(l',r',\chi')$ are equivalent by a $\cb[\partial]$-module homomorphism $\beta:R\rightarrow H$, then we have the following commutative diagram
$$\xymatrix
{0\ar[r]& H\ar[r]\ar@{=}[d]& R\oplus H_{(l,r,\chi)}\ar[r]\ar[d]^{\varphi}& R\ar[r]\ar@{=}[d]& 0\\
0\ar[r]& H\ar[r]& R\oplus H_{(l',r',\chi')}\ar[r]& R\ar[r]& 0,}$$
where $\varphi:R\oplus H_{(l,r,\chi)}\rightarrow R\oplus H_{(l',r',\chi')}$ is defined by
$$\varphi(a,h)=\big(a,h-\beta(a)\big),$$
which implies that the above two extensions are equivalent.
\end{proof}

\bigskip

\section{Extensions in terms of Deligne groupoids}

This section gives the interpretation of $H^2_{nab}(R,H)$, the non-abelian $2$nd cohomology of $R$ with values in $H$, as the Deligne groupoid of a suitable differential graded Lie algebra $\mathfrak{L}$. Therefore, we introduce the differential graded Lie algebra $\mathfrak{L}$ and show that the set of its Maurer-Cartan elements in bijection with the set of extensions of $R$ by $H$.

Given an extension $E$ of $R_1$ by $H_1$, and a section $s:R_1\rightarrow E$ of the extension $E$, we have the following diagram
$$\xymatrix
{0\ar[r]&H_1\ar[r]^\alpha\ar@[r][d]& E\ar[r]^p\ar@[r][d]& R_1\ar[r]\ar@[r][d]& 0\\
0\ar[r]& H\ar[r]^i& R\oplus H\ar[r]^\pi& R\ar[r]& 0,}$$
where $H$ is the image of $H_1$ in $E$ and $R=s(R_1)$ is an arbitrary supplementary in $E$. Then in the category of $\cb[\partial]$-modules, the above diagram is commutative.
Using these isomorphisms, we can transfer the Leibniz conformal algebra structure of $E, R_1$ and $H_1$ to $R\oplus H, R$ and $H$. And the corresponding $\lambda$-bracket on $R\oplus H$ is defined by $\rho$,

We first introduce a notation. Consider the canonical projection
$$P_{A_1,\ldots,A_n}:(R\oplus H)^{\otimes n}\rightarrow A_1\otimes \cdots\otimes A_n,$$
where $A_i\in\{R,H\}$.

Let $\varphi:(R\oplus H)^{\otimes n}\rightarrow (R\oplus H)[\lambda_1,\lambda_2,\ldots,\lambda_{n-1}]$ be a conformal sesquilinear map. We denote
$$\varphi_{A_1,\ldots,A_n}^{A_{n+1}}:=P_{A_{n+1}}\circ \varphi\circ P_{A_1,\ldots,A_n}.$$

Using above notations, all the components of $\lambda$-bracket $\rho$ on $E$ are $\rho^R_{RR},$ $\rho^R_{RH},$ $\rho^R_{HR},$ $\rho^R_{HH},$ $\rho^H_{RR}, $ $\rho^H_{RH},$ $\rho^H_{HR},$ $\rho^H_{HH}$. Then we have

\begin{pro}\label{comp}
$(1)$ The components $\rho^R_{RH},$ $\rho^R_{HR},$ $\rho^R_{HH}$ vanish, $\rho^R_{RR}$ is exactly $[\cdot,\cdot]_R$ and $\rho^H_{HH}$ is $[\cdot,\cdot]_H$ similarly.

$(2)$  Denote by $(l,r,\chi)$ the non-abelian $2$-cocycle on $R_1$ with values in $H_1$ corresponding to the extension $E$. Then the components $\rho^H_{RR},$ $\rho^H_{RH}$, $\rho^H_{HR}$ are $\chi$, $l$, $r$  respectively in fact.
\end{pro}
\begin{proof}
$(1)$ For any $a_1=s(a)\in R, h_1=i(h)\in H$, we have
$$p([{a_1}_\lambda h_1])=p([s(a)_\lambda i(h)])=[p\circ s(a)_\lambda p\circ i(h)]=0,$$
then $[{a_1}_\lambda h_1]\in \rm Ker(p)[\lambda]=\rm Im(i)[\lambda]$, which implies ${(\rho^R_{RH})}_\lambda(R\oplus H)\in H[\lambda].$ Thus $\rho^R_{RH}$ vanishes. Similarly $\rho^R_{HR}$, $\rho^R_{HH}$ vanish.

Since the projection $p$ has inverse $s$, we can identify $\rho^R_{RR}, \rho^H_{HH}$ with the $\lambda$-brackets in $R_1, H_1$ respectively.

$(2)$ If $(l,r,\chi)$ is the non-abelian $2$-cocycle on $R_1$ with values in $H_1$ corresponding to the extension $E$, then the $\lambda$-bracket in $R\oplus H$ is given by
\begin{eqnarray*}
\rho_\lambda((a_1,h_1),(a_2,h_2))=\big([{a_1}_\lambda a_2]_R,l(a_1)_\lambda(h_2)+r(a_2)_{-\partial-\lambda}(h_1)+\chi_\lambda(a_1,a_2)+[{h_1}_\lambda h_2]\big),
\end{eqnarray*}
for any $a_1,a_2\in R, h_1,h_2\in H$. Thus,
\begin{eqnarray*}
{(\rho^H_{RR})}_\lambda((a_1,0),(a_2,0))=P_H\big([{a_1}_\lambda a_2],\chi_\lambda(a_1,a_2)\big)=\chi_\lambda(a_1,a_2),\\
{(\rho^H_{RH})}_\lambda((a_1,0),(0,h_2))=P_H\big(0,l(a_1)_\lambda(h_2)\big)=l(a_1)_\lambda(h_2),\\
{(\rho^H_{HR})}_\lambda((0,h_1),(a_2,0))=P_H\big(0,r(a_2)_\lambda(h_1)\big)=r(a_2)_{-\partial-\lambda}(h_1),
\end{eqnarray*}
which imply the above conclusions.
\end{proof}

Let us now analyze the different components of the Jacobiator of $\rho$. Denote
\begin{eqnarray*}
(J_{RRH}^R)_{\lambda,\mu}(e_1,e_2,e_3):=P_R([{a_1}_\lambda[{a_2}_\mu h_3]]-[{a_2}_\mu[{a_1}_\lambda h_3]]-[[{a_1}_\lambda a_2]_{\lambda+\mu} h_3]),
\end{eqnarray*}
where $e_1:=(a_1,h_1),e_2:=(a_2,h_2),e_3:=(a_3,h_3)\in R\oplus H$. Similarly, we can denote other components of the Jacobiator of $\rho$.
\begin{pro}\label{Jacobi}
Under the conditions above, we have that the components of the Jacobiator of $\rho$ $J_{RRH}^R,$ $J_{RHR}^R,$ $J_{HRR}^R,$ $J_{RHH}^R$, $J_{HRH}^R$, $J_{HHR}^R$, $J_{HHH}^R$ vanish, $J_{RRR}^R$ is Jacobi identity on $R$, and $J_{HHH}^R$ is Jacobi identity on $H$ exactly.
\end{pro}
\begin{proof}
For any $t=(e_1,e_2,e_3)$ with $e_1:=a_1+h_1,e_2:=a_2+h_2,e_3:=a_3+h_3\in R\oplus H,$ we obtain that
\begin{eqnarray*}
(J_{RRH}^R)_{\lambda,\mu}(t)=P_R([{a_1}_\lambda[{a_2}_\mu h_3]]-[{a_2}_\mu[{a_1}_\lambda h_3]]-[[{a_1}_\lambda a_2]_{\lambda+\mu} h_3])=0,
\end{eqnarray*}
Similarly, we have $J_{RHR}^R=J_{HRR}^R=J_{RHH}^R=J_{HRH}^R=J_{HHR}^R=J_{HHH}^R=0$. And
\begin{eqnarray*}
(J_{RRR}^R)_{\lambda,\mu}(t)=[{a_1}_\lambda[{a_2}_\mu a_3]_R]_R-[{a_2}_\mu[{a_1}_\lambda a_3]_R]_R-[[{{a_1}_\lambda a_2]_R}_{\lambda+\mu} a_3]_R=0,\\
(J_{HHH}^H)_{\lambda,\mu}(t)=[{h_1}_\lambda[{h_2}_\mu h_3]_H]_H-[{h_2}_\mu[{h_1}_\lambda h_3]_H]_H-[{[{h_1}_\lambda h_2]_H}_{\lambda+\mu} h_3]_H=0.
\end{eqnarray*}
Then we prove the proposition.
\end{proof}
In addition, we can also compute the following items
\begin{eqnarray*}
(J_{RRR}^H)_{\lambda,\mu}(t)&=&P_H([{a_1}_\lambda[{a_2}_\mu a_3]]-[{a_2}_\mu[{a_1}_\lambda a_3]]-[[{a_1}_\lambda a_2]_{\lambda+\mu} a_3])\\
&=&P_H([{a_1}_\lambda([{a_2}_\mu a_3]_R+\chi_\mu(a_2,a_3))]-[{a_2}_\mu([{a_1}_\lambda a_3]_R+\chi_\mu(a_1,a_3))]\\
&&-[([{a_1}_\lambda a_2]_R+\chi_\mu(a_1,a_2))_{\lambda+\mu} a_3])\\
&=&l(a_1)_\lambda\chi_\mu(a_2,a_3)+\chi_\lambda(a_1,[{a_2}_\mu a_3]_R)-l(a_2)_\mu\chi_\lambda(a_1,a_3)-\chi_\mu(a_2,[{a_1}_\lambda a_3]_R)\\
&&-r(a_3)_{-\partial-\lambda-\mu}\chi_\mu(a_1,a_2)-\chi_{\lambda+\mu}([{a_1}_\lambda a_2]_R,a_3),
\end{eqnarray*}
\begin{eqnarray*}
(J_{RRH}^H)_{\lambda,\mu}(t)&=&P_H([{a_1}_\lambda[{a_2}_\mu h_3]]-[{a_2}_\mu[{a_1}_\lambda h_3]]-[[{a_1}_\lambda a_2]_{\lambda+\mu} h_3])\\
&=&P_H([{a_1}_\lambda (l(a_2)_\mu h_3)]-[{a_2}_\mu (l(a_1)_\lambda h_3)]-[([{a_1}_\lambda a_2]_R+\chi_\lambda(a_1,a_2))_{\lambda+\mu} h_3])\\
&=&l(a_1)_\lambda (l(a_2)_\mu h_3)-l(a_2)_\mu (l(a_1)_\lambda h_3)-l([{a_1}_\lambda a_2]_R)_{\lambda+\mu} h_3+[\chi_\lambda(a_1,a_2)_{\lambda+\mu} h_3],
\end{eqnarray*}
\begin{eqnarray*}
(J_{RHR}^H)_{\lambda,\mu}(t)&=&P_H([{a_1}_\lambda[{h_2}_\mu a_3]]-[{h_2}_\mu[{a_1}_\lambda a_3]]-[[{a_1}_\lambda h_2]_{\lambda+\mu} a_3])\\
&=&P_H([{a_1}_\lambda (r(a_3)_{-\partial-\mu} h_2)]-[{h_2}_\mu ([{a_1}_\lambda a_3]_R+\chi_\lambda(a_1,a_3))]-[(l(a_1)_\lambda h_2)_{\lambda+\mu} a_3])\\
&=&l(a_1)_\lambda(r(a_3)_{-\partial-\mu} h_2)-r([{a_1}_\lambda a_3]_R)_{-\partial-\mu}h_2-[{h_2}_\mu\chi_\lambda(a_1,a_3))]-r(a_3)_{-\partial-\lambda-\mu}(l(a_1)_\lambda h_2),
\end{eqnarray*}
\begin{eqnarray*}
(J_{HRR}^H)_{\lambda,\mu}(t)&=&P_H([{h_1}_\lambda[{a_2}_\mu a_3]]-[{a_2}_\mu[{h_1}_\lambda a_3]]-[[{h_1}_\lambda a_2]_{\lambda+\mu} a_3])\\
&=&P_H([{h_1}_\lambda([{a_2}_\mu a_3]_R+\chi_\mu(a_2,a_3))]-[{a_2}_\mu(r(a_3)_{-\partial-\lambda}h_1)]-[(r(a_2)_{\mu}h_1)_{\lambda+\mu} a_3])\\
&=&r([{a_2}_\mu a_3]_R)_{-\partial-\lambda}h_1+[{h_1}_\lambda\chi_\mu(a_2,a_3)]-l(a_2)_\mu(r(a_3)_{-\partial-\lambda}h_1)-r(a_3)_{-\partial-\lambda-\mu}(r(a_2)_{\mu}h_1),
\end{eqnarray*}
\begin{eqnarray*}
(J_{RHH}^H)_{\lambda,\mu}(t)&=&P_H([{a_1}_\lambda[{h_2}_\mu h_3]]-[{h_2}_\mu[{a_1}_\lambda h_3]]-[[{a_1}_\lambda h_2]_{\lambda+\mu} h_3])\\
&=&l(a_1)_\lambda[{h_2}_\mu h_3]-[{h_2}_\mu (l(a_1)_\lambda h_3)]-[(l(a_1)_\lambda h_2)_{\lambda+\mu} h_3],
\end{eqnarray*}
\begin{eqnarray*}
(J_{HRH}^H)_{\lambda,\mu}(t)&=&P_H([{h_1}_\lambda[{a_2}_\mu h_3]]-[{a_2}_\mu[{h_1}_\lambda h_3]]-[[{h_1}_\lambda a_2]_{\lambda+\mu} h_3])\\
&=&[{h_1}_\lambda (l(a_2)_\mu h_3)]-l(a_2)_\mu[{h_1}_\lambda h_3]-[(r(a_2)_\mu h_1)_{\lambda+\mu} h_3],
\end{eqnarray*}
\begin{eqnarray*}
(J_{HHR}^H)_{\lambda,\mu}(t)&=&P_H([{h_1}_\lambda[{h_2}_\mu a_3]]-[{h_2}_\mu[{h_1}_\lambda a_3]]-[[{h_1}_\lambda h_2]_{\lambda+\mu} a_3])\\
&=&[{h_1}_\lambda (r(a_3)_{-\partial-\mu} h_2)]-[{h_2}_\mu (r(a_3)_{-\partial-\lambda} h_1)]-r(a_3)_{-\partial-\lambda-\mu} [{h_1}_\lambda h_2].
\end{eqnarray*}

Given two Leibniz conformal algebras $R$ and $H$, we get a direct sum $R\oplus H$. It is also a Leibniz conformal algebra under the $\lambda$-bracket
$$[{(a_1,h_1)}_\lambda(a_2,h_2)]_{R\oplus H}=([{a_1}_\lambda a_2]_R,[{h_1}_\lambda h_2]_H).$$
Considering the cohomology complex of $R\oplus H$, we get a differential graded Lie algebra $C(R\oplus H, R\oplus H)$ with the Nijenhuis-Richardson bracket. Also $H$ is an $R\oplus H$-module via the adjoint action of $H$ on itself. Considering the cohomology complex of $R\oplus H$ with coefficients in $A$, we define
$$\mathfrak{L}(R\oplus H, H):=\bigoplus_{i\geq0}\mathfrak{L}^{i}=\bigoplus_{i\geq0}\big(\bigoplus_{m+n=i+1,(m,n)\in{\nb}^*\times \nb}C^{m,n}\big),$$
where $C^{m,n}$ consist of elements of all conformal sesquilinear maps from $(R^{\otimes m}\otimes H^{\otimes n})\oplus (R^{\otimes {m-1}}\otimes H^{\otimes n}\otimes R)\oplus\cdots\oplus (H^{\otimes n}\otimes R^{\otimes m})$ to $H$.

\begin{pro}
With the above notations, the complex $\mathfrak{L}(R\oplus H, H)$ is a sub-differential graded Lie algebra of $C^{\bullet+1}(R\oplus H, R\oplus H)$. We will denote by $(\mathfrak{L}, \delta, [\cdot,\cdot])$ this differential graded Lie algebra. And its degree $0$ part is abelian.
\end{pro}
\begin{proof}
Firstly, for $f\in C^{m,n}, g\in C^{m',n'}$, by the notion of Nijenhuis-Richardson bracket, we have
$[f,g]\in C^{m+m',n+n'-1}.$ Thus, $\mathfrak{L}$ is closed under this bracket. Secondly, for any $f\in C^{m,n}\subseteq \mathfrak{L}^{i}$, we have
$$\delta f=[\rho_R+\rho_H,f]=[\rho_R,f]+[\rho_H,f].$$
Obviously, $\rho_H\in C^{m+1,n},\rho_R\in C^{m,n+1},$ it implies $\delta f\in\mathfrak{L}^{i+1}$.

Finally, note that $\mathfrak{L}^{0}$ is the set of all the $\mathbb{C}[\partial]$-module homomorphisms from $R$ to $H$, we get $\mathfrak{L}^{0}$ is abelian under the Nijenhuis-Richardson bracket.
\end{proof}

Recall that for $\gamma\in \mathfrak{L}^{1}$, $e_1,e_2,e_3\in R\oplus H,$ we have
\begin{eqnarray*}
(\delta\gamma)_{\lambda, \mu}(e_1,e_2,e_3)
&=&[{e_1}_\lambda \gamma_\mu(e_2,e_3)]-[{e_2}_\mu \gamma_\lambda(e_1,e_3)]-[{\gamma_\lambda(e_1,e_2)}_{\lambda+\mu}e_3]\\
&&-\gamma_{\lambda+\mu}([{e_1}_\lambda e_2],e_3)-\gamma_{\mu}(e_2,[{e_1}_\lambda e_3])+\gamma_{\lambda}(e_1,[{e_2}_\mu e_3]),
\end{eqnarray*}
where $\delta$ is the coboundary operator.
\begin{lem}
Let $R_1$ and $H_1$ be Leibniz conformal algebras on $R$ and $H$, then
$$J_{RRR}^H+J_{RRH}^H+J_{RHR}^H+J_{HRR}^H+J_{RHH}^H+J_{HRH}^H+J_{HHR}^H=0\Longleftrightarrow\rho_{RR}^H+\rho_{RH}^H+\rho_{HR}^H=l+r+\chi\in MC(\mathfrak{L}).$$
\end{lem}
\begin{proof}
Let $c=\rho_{RR}^H+\rho_{RH}^H+\rho_{HR}^H\in MC(\mathfrak{L})$. For $e_1:=a_1+h_1,e_2:=a_2+h_2,e_3:=a_3+h_3\in R\oplus H,$ we consider the equality
$${(\delta c+\frac{1}{2}[c,c])}_{\lambda,\mu}(e_1,e_2,e_3)=0.$$

By definition,
\begin{eqnarray*}
(\delta c)_{\lambda, \mu}(e_1,e_2,e_3)
&=&[{e_1}_\lambda c_\mu(e_2,e_3)]-[{e_2}_\mu c_\lambda(e_1,e_3)]-[{c_\lambda(e_1,e_2)}_{\lambda+\mu}e_3]\\
&&-c_{\lambda+\mu}([{e_1}_\lambda e_2],e_3)-c_{\mu}(e_2,[{e_1}_\lambda e_3])+c_{\lambda}(e_1,[{e_2}_\mu e_3]).
\end{eqnarray*}
We will calculate each item one by one in above Formula. Details are as follows:
\begin{eqnarray*}
[{e_1}_\lambda c_\mu(e_2,e_3)]&=&[{h_1}_\lambda\big({(\rho_{RR}^H+\rho_{RH}^H+\rho_{HR}^H)}_\mu(e_2,e_3)\big)]\\
&=&[{h_1}_\lambda\chi_\mu(a_2,a_3)]+[{h_1}_\lambda\big({l(a_2)}_\mu (h_3)\big)]+[{h_1}_\lambda\big({r(a_3)}_{-\partial-\mu} (h_2)\big)].
\end{eqnarray*}
\begin{eqnarray*}
[{e_2}_\mu c_\lambda(e_1,e_3)]&=&[{h_2}_\mu\big({(\rho_{RR}^H+\rho_{RH}^H+\rho_{HR}^H)}_\lambda(e_1,e_3)\big)]\\
&=&[{h_2}_\mu\chi_\lambda(a_1,a_3)]+[{h_2}_\mu\big({l(a_1)}_\lambda (h_3)\big)]+[{h_2}_\mu\big({r(a_3)}_{-\partial-\lambda} (h_1)\big)].
\end{eqnarray*}
\begin{eqnarray*}
[{(c_\lambda(e_1,e_2))}_{\lambda+\mu}e_3]&=&[{\big({(\rho_{RR}^H+\rho_{RH}^H+\rho_{HR}^H)}_\lambda(e_1,e_2)\big)}_{\lambda+\mu}h_3]\\
&=&[{\big(\chi_\lambda(a_1,a_2)\big)}_{\lambda+\mu}h_3]+[{\big({l(a_1)}_\lambda (h_2)\big)}_{\lambda+\mu}h_3]+[{\big({r(a_2)}_{\mu} (h_1)\big)}_{\lambda+\mu}h_3].
\end{eqnarray*}
\begin{eqnarray*}
c_{\lambda+\mu}([{e_1}_\lambda e_2],e_3)&=&{(\rho_{RR}^H+\rho_{RH}^H+\rho_{HR}^H)}_{\lambda+\mu}([{e_1}_\lambda e_2],e_3)\\
&=&\chi_{\lambda+\mu}([{a_1}_\lambda a_2],a_3)+{l([{a_1}_\lambda a_2])}_{\lambda+\mu}(h_3)+{r(a_3)}_{-\partial-\lambda-\mu}([{h_1}_\lambda h_2]).
\end{eqnarray*}
\begin{eqnarray*}
c_{\mu}(e_2,[{e_1}_\lambda e_3])&=&{(\rho_{RR}^H+\rho_{RH}^H+\rho_{HR}^H)}_{\mu}(e_2,[{e_1}_\lambda e_3])\\
&=&\chi_{\mu}(a_2,[{a_1}_\lambda a_3])+{l(a_2)}_{\mu}([{h_1}_\lambda h_3])+{r([{a_1}_\lambda a_3])}_{-\partial-\mu}(h_2).
\end{eqnarray*}
\begin{eqnarray*}
c_{\lambda}(e_1,[{e_2}_\mu e_3])&=&{(\rho_{RR}^H+\rho_{RH}^H+\rho_{HR}^H)}_{\lambda}(e_1,[{e_2}_\mu e_3])\\
&=&\chi_{\lambda}(a_1,[{a_2}_\mu a_3])+{l(a_1)}_{\lambda}([{h_2}_\mu h_3])+{r([{a_2}_\mu a_3])}_{-\partial-\lambda}(h_1).
\end{eqnarray*}
And by the Nijenhuis-Richardson bracket, we have
\begin{eqnarray*}
\frac{1}{2}{[c,c]}_{\lambda,\mu}(e_1,e_2,e_3)&=&{(c\cdot c)}_{\lambda,\mu}(e_1,e_2,e_3)\\
&=&c_{\lambda+\mu}(c_{\lambda}(e_1,e_2),e_3)-c_{\mu}(e_2,c_{\lambda}(e_1,e_3))+c_{\lambda}(e_1,c_{\mu}(e_2,e_3)).
\end{eqnarray*}
And each item is calculated as follows:
\begin{eqnarray*}
&&c_{\lambda+\mu}(c_{\lambda}(e_1,e_2),e_3)\\
&=&{(\rho_{RR}^H+\rho_{RH}^H+\rho_{HR}^H)}_{\lambda+\mu}({(\rho_{RR}^H+\rho_{RH}^H+\rho_{HR}^H)}_{\lambda}(e_1, e_2),e_3)\\
&=&{(\rho_{RR}^H+\rho_{RH}^H+\rho_{HR}^H)}_{\lambda+\mu}(\chi_\lambda(a_1,a_2)+{l(a_1)}_\lambda (h_2)+{r(a_2)}_\mu (h_1),e_3)\\
&=&{r(a_3)}_{-\partial-\lambda-\mu}(\chi_\lambda(a_1,a_2))+{r(a_3)}_{-\partial-\lambda-\mu}({l(a_1)}_\lambda (h_2))+{r(a_3)}_{-\partial-\lambda-\mu}({r(a_2)}_\mu(h_1)).\\
&&c_{\mu}(e_2, c_{\lambda}(e_1,e_3))\\
&=&{(\rho_{RR}^H+\rho_{RH}^H+\rho_{HR}^H)}_{\mu}(e_2,{(\rho_{RR}^H+\rho_{RH}^H+\rho_{HR}^H)}_{\lambda}(e_1,e_3))\\
&=&{(\rho_{RR}^H+\rho_{RH}^H+\rho_{HR}^H)}_{\mu}(e_2,\chi_\lambda(a_1,a_3)+{l(a_1)}_\lambda (h_3)+{r(a_3)}_{-\partial-\lambda} (h_1))\\
&=&{l(a_2)}_{\mu}(\chi_\lambda(a_1,a_3))+{l(a_2)}_{\mu}({l(a_1)}_\lambda (h_3))+{l(a_2)}_{\mu}({r(a_3)}_{-\partial-\lambda}(h_1)).\\
&&c_{\lambda}(e_1,c_{\mu}(e_2,e_3))\\
&=&{(\rho_{RR}^H+\rho_{RH}^H+\rho_{HR}^H)}_{\lambda}(e_1,{(\rho_{RR}^H+\rho_{RH}^H+\rho_{HR}^H)}_{\mu}(e_2,e_3))\\
&=&{(\rho_{RR}^H+\rho_{RH}^H+\rho_{HR}^H)}_{\lambda}(e_1,\chi_\mu(a_2,a_3)+{l(a_2)}_\mu (h_3)+{r(a_3)}_{-\partial-\mu}(h_2))\\
&=&{l(a_1)}_{\lambda}(\chi_\mu(a_2,a_3))+{l(a_1)}_{\lambda}({l(a_2)}_\mu(h_3))+{l(a_1)}_{\lambda}({r(a_3)}_{-\partial-\mu}(h_2)).
\end{eqnarray*}
From above calculations, we can see that
$$(\delta c+\frac{1}{2}[c,c])_{\lambda,\mu}(e_1,e_2,e_3)={(J_{RRR}^H+J_{RRH}^H+J_{RHR}^H+J_{HRR}^H+J_{RHH}^H+J_{HRH}^H+J_{HHR}^H)}_{\lambda,\mu}(e_1,e_2,e_3).$$
Which implies $J_{RRR}^H+J_{RRH}^H+J_{RHR}^H+J_{HRR}^H+J_{RHH}^H+J_{HRH}^H+J_{HHR}^H=0\Longleftrightarrow c\in MC(\mathfrak{L}).$
\end{proof}

\begin{cor} $R,H$ are defined as above. Then
$$Z^2_{nab}(R,H)\cong MC(\mathfrak{L}).$$
\end{cor}
\begin{proof}
For a non-abelian $2$-cocycle $(l,r,\chi)\in Z^2_{nab}(R,H)$, that is $(l,r,\chi)$ satisfies Definition \ref{2-cocycle}, by Proposition \ref{class}, we can define an extension by
$$0 \to H \stackrel{}\to R\oplus H \stackrel{}\to R \to 0.$$
The $\lambda$-bracket $\rho$ defined on $R\oplus H$ is
\begin{eqnarray}
[(a,u)_\lambda(b,v)]=\big([a_\lambda b],l(a)_\lambda(v)+r(b)_{-\partial-\lambda}(u)+\chi_\lambda(a,b)+[u_\lambda v]\big).\label{non-abelian}
\end{eqnarray}
From Propositions \ref{comp} and \ref{Jacobi}, we can obtain the Jacobiator of $\rho$ are $0$, that is $J_{RRR}^H+J_{RRH}^H+J_{RHR}^H+J_{HRR}^H+J_{RHH}^H+J_{HRH}^H+J_{HHR}^H=0$, which means that $\rho_{RR}^H+\rho_{RH}^H+\rho_{HR}^H\in MC(\mathfrak{L})$.
\end{proof}

\begin{thm}
Equivalence in $Z^2_{nab}(R,H)$ coincide with equivalence in $MC(\mathfrak{L})$.
\end{thm}
\begin{proof}
Recall that $\beta$ and $\beta'$ in $MC(\mathfrak{L})$ are equivalent if there exists an $\alpha\in L^0$ such that
$$\beta'=e^{\rm ad\alpha}(\beta)-\frac{e^{\rm ad\alpha}-1}{\rm ad\alpha}\delta(\alpha).$$

Let $\beta:=l+r+\chi$, and $e_1:=a_1+h_1,e_2:=a_2+h_2\in R\oplus H$.
\begin{eqnarray*}
&&{\big(e^{\rm ad\alpha}(l+r+\chi)\big)}_\lambda(e_1,e_2)\\&=&{\big(l+r+\chi+[\alpha,l+r+\chi]\big)}_\lambda(e_1,e_2)\\
&=&l(a_1)_\lambda(h_2)+r(a_2)_{-\partial-\lambda}(h_1)+\chi_\lambda(a_1,a_2)+l(a_1)_\lambda(\alpha(a_2))+r(a_2)_{-\partial-\lambda}(\alpha(a_1)).
\end{eqnarray*}
We now compute that
\begin{eqnarray*}
-\big(\frac{e^{\rm ad\alpha}-1}{\rm ad\alpha}\delta(\alpha)\big)_\lambda(e_1,e_2)&=&-\sum_{n=0}\frac{1}{(n+1)!}{\big((\rm ad\alpha)^n\delta(\alpha)\big)}_\lambda(e_1,e_2)\\
&=&-\frac{1}{2}{\big((\rm ad\alpha)\delta(\alpha)\big)}_\lambda(e_1,e_2)\\
&=&-\alpha([{a_1}_\lambda a_2])+[{\alpha(a_1)}_\lambda h_2]+[{h_1}_\lambda\alpha(a_2)]+[{\alpha(a_1)}_\lambda\alpha(a_2)].
\end{eqnarray*}
Combining these results gives
\begin{eqnarray*}
\beta'_\lambda(e_1,e_2)&=&\beta_\lambda(e_1,e_2)+l(a_1)_\lambda(\alpha(a_2))+r(a_2)_{-\partial-\lambda}(\alpha(a_1))\\
&&+{{\rm ad}^L_{\alpha(a_1)}}_\lambda(h_2)+{{\rm ad}^R_{\alpha(a_2)}}_{-\partial-\lambda}(h_1)+[{\alpha(a_1)}_\lambda\alpha(a_2)]-\alpha([{a_1}_\lambda a_2]).
\end{eqnarray*}

On the other hand, the two cocycles $(l,r,\chi)$ and $(l',r',\chi')$ in $Z^2_{nab}(R,H)$ are equivalent if and only if Eqs. (\ref{cocycle1}) and (\ref{cocycle2}) are satisfied, that is
\begin{eqnarray*}
\beta'_\lambda(e_1,e_2)&=&\beta_\lambda(e_1,e_2)+l(a_1)_\lambda(\alpha(a_2))+r(a_2)_{-\partial-\lambda}(\alpha(a_1))\\
&&+{{\rm ad}^L_{\alpha(a_1)}}_\lambda(h_2)+{{\rm ad}^R_{\alpha(a_2)}}_{-\partial-\lambda}(h_1)+[{\alpha(a_1)}_\lambda\alpha(a_2)]-\alpha([{a_1}_\lambda a_2]),
\end{eqnarray*}
where $\beta'=l'+r'+\chi'$ and $\beta=l+r+\chi$. Hence the theorem is proved.
\end{proof}

\bigskip

\section{The inducibility problem of automorphisms}

Let $R$ and $H$ be two Leibniz conformal algebras, $\mathcal{E}:0 \to H \stackrel{\alpha}\to E \stackrel{p}\to R \to 0$ be a non-abelian extension of $R$ by $H$ with a section $s$. Denote the corresponding non-abelian $2$-cocycle of $\mathcal{E}$ by $(r,l,\chi)$. From the Section $4$, we can identify $E$ with $R\oplus H$ under the $\lambda$-bracket defined as Eq. (\ref{non-abelian}), identify $\alpha$ and $p$ with an injection and a projection respectively. Denote
$$\Aut_H(E)=\{f\in\Aut(E)\mid f(H)=H\}.$$

Obviously, if $f\in \Aut_H(E)$, then $f|_H\in\Aut(H)$. For any $f\in \Aut_H(E)$, define a $\cb[\partial]$-module homomorphism
$$\hat{f}=p\circ f\circ s:R\rightarrow R.$$
Firstly, $\hat{f}$ is independent of the choice of $s$. Since if $s$ and $s'$ are two sections of the extension $\mathcal{E}$, for any $a\in R$, we have $p\circ(s-s')(a)=0$, that is
$(s-s')(a)\in \Ker(p)=\rm Im(\alpha)$. Then there is a $h\in H$ such that $p\circ f\circ (s-s')(a)=p\circ f\circ\alpha(h)=0.$ Similarly, we can obtain that $\hat{f}$ is injective and it is a surjection obviously. Hence $\hat{f}$ is a bijection on $R$. Furthermore, for $a,b\in R$,
\begin{eqnarray*}
\hat{f}([a_\lambda b])&=&(p\circ f\circ s)([a_\lambda b])=(p\circ f)\big([{s(a)}_\lambda s(b)]-\chi_\lambda(a,b)\big)\\
&=&(p\circ f)\big([{s(a)}_\lambda s(b)]\big)=[{\hat{f}(a)}_\lambda \hat{f}(b)].
\end{eqnarray*}
This means $\hat{f}\in\Aut(R).$ Then we obtain a group homomorphism
$$\kappa:\Aut_H(E)\rightarrow\Aut(H)\times\Aut(R)\ \ \ \ \ \ f\mapsto (f|_H,\hat{f}).$$

\begin{defi}
Let $R$ and $H$ be two Leibniz conformal algebras, $\mathcal{E}:0 \to H \stackrel{\alpha}\to E \stackrel{p}\to R \to 0$ be a non-abelian extension of $R$ by $H$ with a section $s$, and $(r,l,\chi)$ be the corresponding non-abelian $2$-cocycle. A pair $(\gamma,\varphi)\in\Aut(H)\times\Aut(R)$ is said to be {\rm inducible} if there exists a map $f\in\Aut_H(E)$ such that $(\gamma,\varphi)=(f|_H,\hat{f})$.
\end{defi}

For a pair $(\gamma,\varphi)\in\Aut(H)\times\Aut(R)$ and a non-abelian $2$-cocycle $(r,l,\chi)$, we can get a new pair $(l^{\gamma,\varphi}, r^{\gamma,\varphi},\chi^{\gamma,\varphi})$ details as follows:
\begin{eqnarray*}
l^{\gamma,\varphi}:R\rightarrow {\rm Cend}(H),\ \ \ \ \ l^{\gamma,\varphi}(a)_\lambda(h)=\gamma\big(l(\varphi^{-1}(a))_\lambda \gamma^{-1}(h)\big),\\
r^{\gamma,\varphi}:R\rightarrow {\rm Cend}(H),\ \ \ \ \ r^{\gamma,\varphi}(a)_\lambda(h)=\gamma\big(r(\varphi^{-1}(a))_\lambda \gamma^{-1}(h)\big),\\
\chi^{\gamma,\varphi}:R\times R\rightarrow H[\lambda],\ \ \ \ \chi^{\gamma,\varphi}_\lambda(a,b)=\gamma\circ\chi_\lambda\big(\varphi^{-1}(a),\varphi^{-1}(b)\big),
\end{eqnarray*}
for any $a,b\in R, h\in H.$ Then we can check that $(l^{\gamma,\varphi},r^{\gamma,\varphi},\chi^{\gamma,\varphi})$ is also a non-abelian $2$-cocycle.

\begin{thm}
Let $R,H,\mathcal{E},(r,l,\chi),(l^{\gamma,\varphi},r^{\gamma,\varphi},\chi^{\gamma,\varphi})$ be defined as above. Then a pair $(\gamma,\varphi)\in\Aut(H)\times\Aut(R)$ is inducible if and only if the non-abelian $2$-cocycles $(r,l,\chi)$ and $(l^{\gamma,\varphi},r^{\gamma,\varphi},\chi^{\gamma,\varphi})$ are equivalent.
\end{thm}
\begin{proof}
If the pair $(\gamma,\varphi)$ is inducible, there exists a map $f\in\Aut_H(E)$ such that $f|_H=\gamma$ and $\hat{f}=p\circ f\circ s=\varphi$, where $s$ is a section of the extension $\mathcal{E}$. Obviously, for any $a\in R$, we have $p\big(s(a)-f\circ s\circ\varphi^{-1}(a)\big)=a-a=0$, which implies $s(a)-f\circ s\circ\varphi^{-1}(a)\in {\rm Ker}(p)=\rm Im(\alpha)\cong H$. Thus we can define a $\cb[\partial]$-module homomorphism
$$\beta:R\rightarrow H\ \ \ \ a\mapsto s(a)-f\circ s\circ\varphi^{-1}(a).$$
Then by Eq. (\ref{extension}) for any $h\in H$ and $a_1,a_2\in R$, we compute that
\begin{eqnarray*}
\big(l-l^{\gamma,\varphi}\big)(a)_\lambda(h)=[s(a)_\lambda h]-\gamma\big(l(\varphi^{-1}(a))_\lambda\gamma^{-1}(h)\big)=[s(a)_\lambda h]-[\gamma\circ s\circ\varphi^{-1}(a)_\lambda h]=[\beta(a)_\lambda h],
\end{eqnarray*}
\begin{eqnarray*}
\big(r-r^{\gamma,\varphi}\big)(a)_\lambda(h)=[h_{-\partial-\lambda}s(a)]-[h_{-\partial-\lambda}\gamma\circ s\circ\varphi^{-1}(a)]=[h_{-\partial-\lambda}\beta(a)],
\end{eqnarray*}
and
\begin{eqnarray*}
&&\big(\chi-\chi^{\gamma,\varphi}\big)_\lambda(a_1,a_2)\\
&=&\chi_\lambda(a_1,a_2)-\gamma\circ\chi_\lambda\big(\varphi^{-1}(a_1),\varphi^{-1}(a_2)\big)\\
&=&[{s(a_1)}_\lambda s(a_2)]-s([{a_1}_\lambda a_2])-\gamma\big([s(\varphi^{-1}(a_1))_\lambda s(\varphi^{-1}(a_2))]-s([\varphi^{-1}(a_1)_\lambda\varphi^{-1}(a_2)])]\big)\\
&=&[{s(a_1)}_\lambda s(a_2)]-s([{a_1}_\lambda a_2])-f\big([s(\varphi^{-1}(a_1))_\lambda s(\varphi^{-1}(a_2))]-s([\varphi^{-1}(a_1)_\lambda\varphi^{-1}(a_2)]\big)\\
&=&[\big(f\circ s\circ\varphi^{-1}(a_1)\big)_\lambda \big(s(a_2)-f\circ s\circ\varphi^{-1}(a_2)\big)]+[{\big(s(a_1)-f\circ s\circ\varphi^{-1}(a_1)\big)}_\lambda \big(f\circ s\circ\varphi^{-1}(a_2)]\\
&&+[{\big(s(a_1)-f\circ s\circ\varphi^{-1}(a_1)\big)}_\lambda \big(s(a_2)-f\circ s\circ\varphi^{-1}(a_2)\big)]-(s-f\circ s\circ\varphi^{-1})([{a_1}_\lambda a_2])\\
&=&[\big(f\circ s\circ\varphi^{-1}(a_1)\big)_\lambda \beta(a_2)]+[\beta(a_1)_\lambda \big(f\circ s\circ\varphi^{-1}(a_2)]
+[\beta(a_1)_\lambda \beta(a_2)]-\beta([{a_1}_\lambda a_2]),
\end{eqnarray*}
where $f\circ s\circ\varphi^{-1}$ is also a section of the extension $\mathcal{E}$, then we have that the above equation is equal to
$$l(a_1)_\lambda \beta(a_2)+r(a_2)_{-\partial-\lambda}\beta(a_1)+[\beta(a_1)_\lambda \beta(a_2)]-\beta([{a_1}_\lambda a_2]).$$
Hence $(r,l,\chi)$ and $(l^{\gamma,\varphi},r^{\gamma,\varphi},\chi^{\gamma,\varphi})$ are equivalent.

Conversely, if the non-abelian $2$-cocycles $(r,l,\chi)$ and $(l^{\gamma,\varphi},r^{\gamma,\varphi},\chi^{\gamma,\varphi})$ are equivalent by a $\cb[\partial]$-module homomorphism $\beta:R\rightarrow H$ satisfying Eqs. (\ref{cocycle1}) and (\ref{cocycle2}), then we have the following commutative diagram
$$\xymatrix
{0\ar[r]& H\ar[r]\ar@{=}[d]& R\oplus H_{(r,l,\chi)}\ar[r]\ar[d]^{\theta}& R\ar[r]\ar@{=}[d]& 0\\
0\ar[r]& H\ar[r]& R\oplus H_{(l^{\gamma,\varphi},r^{\gamma,\varphi},\chi^{\gamma,\varphi})}\ar[r]& R\ar[r]& 0,}$$
where $\theta:R\oplus H_{(r,l,\chi)}\rightarrow R\oplus H_{(l^{\gamma,\varphi},r^{\gamma,\varphi},\chi^{\gamma,\varphi})}$ is defined by
$$\theta(a,h)=\big(a,h-\beta(a)\big).$$
Since $E\cong R\oplus H_{(r,l,\chi)}$ as Leibniz conformal algebras, we write each element as $(a,h)$ in $E$. Define a $\cb[\partial]$-module homomorphism
$$f:E\rightarrow E\ \ \ \ \ \ (a,h)\mapsto\theta(\varphi(a),\gamma(h))=(\varphi(a),\gamma(h)-\beta(\varphi(a))),$$
for any $a\in R,h\in H$. Then $f$ is a bijection on $E$. Furthermore, for any $(a_1,h_1),(a_2,h_2)\in E$, we have
\begin{eqnarray*}
&&[f(a_1,h_1)_\lambda f(a_2,h_2)]\\
&=&\theta\big([(\varphi(a_1),\gamma(h_1))_\lambda(\varphi(a_2),\gamma(h_2))]\big)\\
&=&\theta\big([\varphi(a_1)_\lambda\varphi(a_2)],l^{\gamma,\varphi}(\varphi(a_1))_\lambda\gamma(h_2)+r^{\gamma,\varphi}(\varphi(a_2))_{-\partial-\lambda}\gamma(h_1)
+\chi^{\gamma,\varphi}_\lambda(\varphi(a_1),\varphi(a_2))+[\gamma(h_1)_\lambda\gamma(h_2)]\big)\\
&=&\theta\big(\varphi([{a_1}_\lambda a_2]),\gamma(l(a_1)_\lambda h_2)+\gamma(r(a_2)_{-\partial-\lambda}h_1)+\gamma(\chi_\lambda(a_1,a_2))+\gamma([{h_1}_\lambda h_2)]\big)\\
&=&f\big([{a_1}_\lambda a_2],l(a_1)_\lambda h_2+r(a_2)_{-\partial-\lambda}h_1+\chi_\lambda(a_1,a_2)+[{h_1}_\lambda h_2]\big)\\
&=&f([(a_1,h_1)_\lambda(a_2,h_2)]),
\end{eqnarray*}
which implies $f\in \Aut(E).$ Obviously, $f|_H=\gamma$ and $\hat{f}=p\circ f\circ s=\varphi.$ So the pair $(\gamma,\varphi)$ is inducible.
\end{proof}

Next, we will give a map which called Wells map to interpret the above theorem. Let $R$ and $H$ be two Leibniz conformal algebras, $\mathcal{E}:0 \to H \stackrel{\alpha}\to E \stackrel{p}\to R \to 0$ be a non-abelian extension of $R$ by $H$, and $(l,r,\chi)$ be the corresponding non-abelian $2$-cocycle. Define a map
$$\mathcal{W}:\Aut(H)\times\Aut(R)\rightarrow H^2_{nab}(R,H), \ \ \ \ \ (\gamma,\varphi)\mapsto[(l^{\gamma,\varphi},r^{\gamma,\varphi},\chi^{\gamma,\varphi})-(l,r,\chi)].$$
With the above definition, it is easy to obtain the following conclusion.

\begin{cor}\label{inducible}
Let $\mathcal{E}:0 \to H \stackrel{\alpha}\to E \stackrel{p}\to R \to 0$ be a non-abelian extension of $R$ by $H$. Then a pair $(\gamma,\varphi)\in\Aut(H)\times\Aut(R)$ is inducible if and only if $\mathcal{W}(\gamma,\varphi)=0$.
\end{cor}

Let $\mathcal{E}:0 \to H \stackrel{\alpha}\to E \stackrel{p}\to R \to 0$ be a non-abelian extension of $R$ by $H$. Define
$$\Aut^{\rm id}_H(E):=\{f\in\Aut_H(E)\mid\kappa(f)=({\rm id}_H,{\rm id}_R)\}.$$
Then we have the following corollary.
\begin{cor}
The sequence of groups
$$1 \to \Aut^{\rm id}_H(E) \stackrel{\iota}\to \Aut_H(E) \stackrel{\kappa}\to \Aut(H)\times\Aut(R) \stackrel{\mathcal{W}}\to H^2_{nab}(R,H)$$
is an exact sequence, where $\iota$ is the inclusion map.
\end{cor}
\begin{proof}
By the definition of $\Aut^{\rm id}_H(E)$ and Corollary \ref{inducible}, it is easy to prove the sequence is exact at $\Aut_H(E)$ and $\Aut(H)\times\Aut(R)$.
\end{proof}
Moreover, we define that the following two sets
$$\Aut^{H}_H(E):=\{f\in\Aut_H(E)\mid f|_H={\rm id}_H\},$$
$$\Aut^{R}_H(E):=\{f\in\Aut_H(E)\mid \hat{f}={\rm id}_R\}.$$
Similar to the above corollary, we have

\begin{pro}
Let $\mathcal{E}:0 \to H \stackrel{\alpha}\to E \stackrel{p}\to R \to 0$ be a non-abelian extension of $R$ by $H$. Then the following two sequences of groups
$$1 \to \Aut^{\rm id}_H(E) \stackrel{\iota_H}\to \Aut^R_H(E) \stackrel{\kappa_H}\to \Aut(H) \stackrel{\mathcal{W}_H}\to H^2_{nab}(R,H),$$
$$1 \to \Aut^{\rm id}_H(E) \stackrel{\iota_R}\to \Aut^H_H(E) \stackrel{\kappa_R}\to \Aut(R) \stackrel{\mathcal{W}_R}\to H^2_{nab}(R,H),$$
are exact.
\end{pro}

Net we recall the abelian extension. Let $R$ be a Leibniz conformal algebra and $H$ be a $\cb[\partial]$-module which can be viewed as a trivial Leibniz conformal algebra. Given an abelian extension $\mathcal{E}:0 \to H \stackrel{\alpha}\to E \stackrel{p}\to R \to 0$, denote by $\Aut(H)$ the set of all $\cb[\partial]$-module automorphisms and $\Aut_\bullet(H,R)$ the set of all pair $(\gamma,\varphi)\in \Aut(H)\times\Aut(R)$ satisfying
$$\gamma(a\triangleright_\lambda h)=\varphi(a)\triangleright_\lambda \gamma(h),\ \ \ \  \gamma(h\triangleleft_\lambda a)=\gamma(h)\triangleleft_\lambda \varphi(a)$$
for any $a\in R,h\in H$, where $(\triangleleft, \triangleright)$ is the $R$-module structure on $H$ induced by the abelian extension $\mathcal{E}$, i.e. $a\triangleright_\lambda h=[s(a)_\lambda h], h\triangleleft_\lambda a=[h_\lambda s(a)]$ with a section $s$. And the $2$-cocycle $\chi:R\times R\rightarrow H[\lambda]$ induced by $\mathcal{E}$ is $\chi_\lambda(a,b)=[s(a)_\lambda s(b)]-s([a_\lambda b])$. Define $\chi^{\gamma,\varphi}:R\times R\rightarrow H[\lambda]$ as follows:
$$\chi^{\gamma,\varphi}_\lambda(a,b)=\gamma\circ\chi_\lambda\big(\varphi^{-1}(a),\varphi^{-1}(b)\big),$$
for any $a,b\in R.$ Furthermore, if $(\gamma,\varphi)\in\Aut_\bullet(H,R)$, for any $a,b,c\in R$. we have
$$\delta(\chi^{\gamma,\varphi})_{\lambda,\mu}(a,b,c)=\gamma\circ\delta(\chi)_{\lambda,\mu}(\varphi^{-1}(a),\varphi^{-1}(b),\varphi^{-1}(c))=0.$$
That is $\chi^{\gamma,\varphi}\in Z^2(R,H)$. Then we have the following theorem.

\begin{thm}
Let $\mathcal{E}:0 \to H \stackrel{\alpha}\to E \stackrel{p}\to R \to 0$ be an abelian extension of $R$ by $H$. Then a pair $(\gamma,\varphi)\in\Aut_\bullet(H,R)$ is inducible if and only if the $2$-cocycles $\chi$ and $\chi^{\gamma,\varphi}$ are equivalent, if and only if there exists a $\cb[\partial]$-module homomorphism $\beta:R\rightarrow H$ satisfying
\begin{eqnarray*}
\gamma\circ\chi_\lambda(a,b)=\chi_\lambda(\varphi(a),\varphi(b))+{\varphi(a)}\triangleright_\lambda\beta(\varphi(b))
+\beta(\varphi(a))\triangleleft_\lambda\varphi(b)-\beta([\varphi(a)_\lambda \varphi(b)]).
\end{eqnarray*}
\end{thm}

Now define the Wells map for the abelian extensions by
$$\mathcal{\bar{W}}:\Aut_\bullet(H,R)\rightarrow H^2(R,H),\ \ \ \ \ (\gamma,\varphi)\mapsto[\chi^{\gamma,\varphi}-\chi].$$
Let $s$ be a section of an abelian extension $\mathcal{E}:0 \to H \stackrel{\alpha}\to E \stackrel{p}\to R \to 0$. For any $f\in\Aut_H(E)$, we define $\hat{f}=p\circ f\circ s.$
Then we obtain a group homomorphism
$$\bar{\kappa}:\Aut_H(E)\rightarrow\Aut_\bullet(H,R)\ \ \ \ \ \ f\mapsto (f|_H,\hat{f}).$$
Then we have the following fundamental sequence of Wells.
\begin{cor}
The sequence of groups
$$1 \to \Aut^{\rm id}_H(E) \stackrel{\iota}\to \Aut_H(E) \stackrel{\bar{\kappa}}\to \Aut_\bullet(H,R) \stackrel{\mathcal{\bar{W}}}\to H^2(R,H)$$
is exact, where $\iota$ is the inclusion map.
\end{cor}

\bigskip

\section{The extensibility problem of derivations}

Let $R$ be a Leibniz conformal algebra and $H$ be a $\cb[\partial]$-module. Given an abelian extension $\mathcal{E}:0 \to H \stackrel{\alpha}\to E \stackrel{p}\to R \to 0$, there is an $R$-module on $H$ given by $a \triangleright_\lambda h=[s(a)_\lambda h], h \triangleleft_\lambda a=[h_\lambda s(a)]$ with a section $s$ and this action does not depend on the choice of $s$ for $a\in R,h\in H$. Then there is a Leibniz conformal algebra structure on $\cb[\partial]$-module $R\oplus H$, which is given by
$$[(a_1,h_1)_\lambda(a_2,h_2)]=\big([{a_1}_\lambda a_2],{a_1}\triangleright_\lambda h_2+h_1\triangleleft_\lambda a_2),$$
for any $a_1,a_2\in R, h_1,h_2\in H$. Denote this Leibniz conformal algebra by $R\ltimes H$.
Recall that the set of left derivations ${\rm Der}^L(R)$ in Definition \ref{derivation}, and it is a Lie algebra. We view $H$ as a trivial Leibniz conformal algebra. Then ${\rm Der}^L(H)$ just is the set of all $\cb[\partial]$-module homomorphisms of $H$. Then by direct calculations, we have the following lemma.
\begin{lem}
Let $R$ be a Leibniz conformal algebra, $H$ be a $\cb[\partial]$-module and $(d_R,d_H)\in{\rm Der}^L(R)\times{\rm Der}^L(H)$. Then $(d_R,d_H)\in{\rm Der}^L(R\ltimes H)$ if and only if
\begin{eqnarray}
d_H(a\triangleright_\lambda h)=a\triangleright_\lambda d_H(h)+d_R(a)\triangleright_\lambda h,\label{Der}\\
d_H(h\triangleleft_\lambda a)=h\triangleleft_\lambda d_R(a)+d_H(h)\triangleleft_\lambda a,\label{Der1}
\end{eqnarray}
for $ a\in R,h\in H.$
\end{lem}
\begin{proof}
It is easy to prove this lemma by the definition of derivation.
\end{proof}
Denote by
\begin{eqnarray*}
{\rm g}(R,H):=\{(d_R,d_H)\in{\rm Der}^L(R)\times{\rm Der}^L(H)\mid(d_R,d_H) \ \ satisfies \ Eqs. \ (\ref{Der}), (\ref{Der1})\}.
\end{eqnarray*}
Then we can check that ${\rm g}(R,H)$ is a Lie subalgebra of ${\rm Der}^L(R\ltimes H)$. The following we can define an action of ${\rm g}(R,H)$ on the space of $2$-cochains $C^2(R,H)$ by
$$\Phi(d_R,d_H)(f)_\lambda(a,b)=d_H(f_\lambda(a,b))-f_\lambda(d_R(a),b)-f_\lambda(a,d_R(b)),$$
for any $a,b \in R, (d_R,d_H)\in {\rm g}(R,H), f\in C^2(R,H)$. Then we have the following proposition.

\begin{pro}
With above notations, we have $$\delta(\Phi(d_R,d_H)(f))_{\lambda,\mu}(a,b,c)=0,$$
for any $a,b,c\in R, f\in Z^2(R,H)$.
\end{pro}
\begin{proof}
For any $f\in Z^2(R,H)$, i.e. $\delta(f)_{\lambda,\mu}(a,b,c)=0$, that is $f$ satisfies
\begin{eqnarray*}
&&a\triangleright_\lambda f_\mu(b,c)-b\triangleright_\mu f_\lambda(a,c)-f_\lambda(a,b)\triangleleft_{\lambda+\mu} c\\
&&\ \ \ \ \ \ \ \ \ \ \ \ \ \ \ -f_{\lambda+\mu}([a_\lambda b],c)-f_{\mu}(b,[a_\lambda c])+f_{\lambda}(a,[b_\mu c])=0,
\end{eqnarray*}
for any $a,b,c\in R.$ By direct calculation, we have
\begin{eqnarray*}
&&\delta(\Phi(d_R,d_H)(f))_{\lambda,\mu}(a,b,c)\\
&=&a\triangleright_\lambda (\Phi(d_R,d_H)(f))_\mu(b,c)-b\triangleright_\mu (\Phi(d_R,d_H)(f))_\lambda(a,c)-(\Phi(d_R,d_H)(f))_\lambda(a,b)\triangleleft _{\lambda+\mu}c\\
&&-(\Phi(d_R,d_H)(f))_{\lambda+\mu}([a_\lambda b],c)-(\Phi(d_R,d_H)(f))_{\mu}(b,[a_\lambda c])+(\Phi(d_R,d_H)(f))_{\lambda}(a,[b_\mu c])\\
&=&a\triangleright_\lambda d_H(f_\mu(b,c))-a\triangleright_\lambda f_\mu(d_R(b),c)-a\triangleright_\lambda f_\mu(b,d_R(c))\\
&&-b\triangleright_\mu d_H(f_\lambda(a,c))+b\triangleright_\mu f_\lambda(d_R(a),c)+b\triangleright_\mu f_\lambda(a,d_R(c))\\
&&-d_H(f_\lambda(a,b))\triangleleft_{\lambda+\mu}c+f_\lambda(d_R(a),b)\triangleleft_{\lambda+\mu}c+f_\lambda(a,d_R(b))\triangleleft_{\lambda+\mu}c\\
&&-d_H(f_{\lambda+\mu}([a_\lambda b],c))+f_{\lambda+\mu}(d_R([a_\lambda b]),c)+f_{\lambda+\mu}([a_\lambda b],d_R(c))\\
&&-d_H(f_{\mu}(b,[a_\lambda c]))+f_{\mu}(d_R(b),[a_\lambda c])+f_{\mu}(b,d_R([a_\lambda c]))\\
&&+d_H(f_{\lambda}(a,[b_\mu c]))-f_{\lambda}(d_R(a),[b_\mu c])-f_{\lambda}(a,d_R([b_\mu c])).
\end{eqnarray*}
Since $d_R$ is a left derivation of $R$ and $(d_R,d_H)$ satisfies Eqs. (\ref{Der}) and (\ref{Der1}), then we obtain
\begin{eqnarray}
f_{\lambda+\mu}(d_R([a_\lambda b]),c)=f_{\lambda+\mu}([d_R(a)_\lambda b],c)+f_{\lambda+\mu}([a_\lambda d_R(b)],c),\label{dH}
\end{eqnarray}
\begin{eqnarray}
a\triangleright_\lambda d_H(f_\mu(b,c))= d_H(a\triangleright_\lambda f_\mu(b,c))-d_R(a)\triangleright_\lambda f_\mu(b,c),\label{dR}\\
d_H(f_\lambda(a,b))\triangleleft_{\lambda+\mu}c=d_H(f_\lambda(a,b)\triangleleft_{\lambda+\mu}c)-f_\lambda(a,b)\triangleleft_{\lambda+\mu}d_R(c).\label{dR1}
\end{eqnarray}
Combining with Eqs. (\ref{dH}), (\ref{dR}) and (\ref{dR1}), it is easy to prove $\delta(\Phi(d_R,d_H)(f))_{\lambda,\mu}(a,b,c)=0,$ which implies the proposition.
\end{proof}
In addition, we can also obtain that for any $f\in C^1(R,H)$, $\Phi(d_R,d_H)(\delta(f))=\delta(d_H\circ f-f\circ d_R).$
Then $\Phi$ can induce a linear map still denoted as $\Phi:{\rm g}(R,H)\rightarrow gl(H^2(R,H))$ as follows.
$$\Phi(d_R,d_H)([f])=[\Phi(d_R,d_H)(f)],$$
for any $(d_R,d_H)\in{\rm g}(R,H)$ and $[f]\in H^2(R,H)$. Furthermore, $\Phi$ is a Lie homomorphism in fact.

\begin{pro}
The map $\Phi$ defines a representation of Lie algebra ${\rm g}(R,H)$ on the space $H^2(R,H)$.
\end{pro}
\begin{proof}
For any $(d_R,d_H), (d'_R,d'_H)\in{\rm g}(R,H)$ and $a,b\in R, [f]\in H^2(R,H)$, we have
\begin{eqnarray*}
&&\big(\Phi(d_R,d_H)\circ\Phi(d'_R,d'_H)([f])\big)_\lambda(a,b)\\
&=&(\Phi(d_R,d_H)([\Phi(d'_R,d'_H)(f)]))_\lambda(a,b)\\
&=&[\Phi(d_R,d_H)\circ\Phi(d'_R,d'_H)(f)]_\lambda(a,b)\\
&=&d_H\big(\Phi(d'_R,d'_H)(f)_\lambda(a,b)\big)-\Phi(d'_R,d'_H)(f)_\lambda(d_R(a),b)-\Phi(d'_R,d'_H)(f)_\lambda(a,d_R(b))\\
&=&d_H\circ d'_H(f_\lambda(a,b))-d_H(f_\lambda(d'_R(a),b))-d_H(f_\lambda(a,d'_R(b)))\\
&&-d'_H(f_\lambda(d_R(a),b))+f_\lambda(d'_R\circ d_R(a),b)+f_\lambda(d_R(a),d'_R(b))\\
&&-d'_H(f_\lambda(a,d_R(b)))+f_\lambda(d'_R(a),d_R(b))+f_\lambda(a,d'_R\circ d_R(b)).
\end{eqnarray*}
Similarly,
\begin{eqnarray*}
&&(\Phi(d'_R,d'_H)\circ\Phi(d_R,d_H)([f]))_\lambda(a,b)\\
&&=d'_H\circ d_H(f_\lambda(a,b))-d'_H(f_\lambda(d_R(a),b))-d'_H(f_\lambda(a,d_R(b)))-d_H(f_\lambda(d'_R(a),b))\\
&&\ \ +f_\lambda(d_R\circ d'_R(a),b)+f_\lambda(d'_R(a),d_R(b))-d_H(f_\lambda(a,d'_R(b)))+f_\lambda(d_R(a),d'_R(b))+f_\lambda(a,d_R\circ d'_R(b)).
\end{eqnarray*}
Then we get
\begin{eqnarray*}
&&(\Phi(d_R,d_H)\circ\Phi(d'_R,d'_H)([f])-\Phi(d'_R,d'_H)\circ\Phi(d_R,d_H)([f]))_\lambda(a,b)\\
&=&[d_H,d'_H](f_\lambda(a,b))-f_\lambda([d_R,d'_R](a),b)-f_\lambda(a,[d_R,d'_R](b))\\
&=&\Phi\big([(d_R,d_H),(d'_R,d'_H)]\big)([f])_\lambda(a,b).
\end{eqnarray*}
Hence $\Phi$ is a Lie algebra homomorphism.
\end{proof}

Now we consider the extensibility problem of derivations about an abelian extension of Leibniz conformal algebras.

\begin{defi}
Let $\mathcal{E}:0 \to H \stackrel{\alpha}\to E \stackrel{p}\to R \to 0$ be an abelian extension of Leibniz conformal algebra $R$ by a $\cb[\partial]$-module $H$. A pair $(d_R,d_H)\in{\rm Der}(R)\times{\rm Der}(H)$ is called an {\rm  extensible} pair if there is a derivation $d_E\in{\rm Der}(E)$ such that $d_E\circ \alpha=\alpha\circ d_H$ and $d_R\circ p=p\circ d_E$.
\end{defi}

Recall that $\mathcal{E}:0 \to H \stackrel{\alpha}\to E \stackrel{p}\to R \to 0$ is an abelian extension of $R$ by $H$. Let $s$ be a section of $\mathcal{E}$ and $\chi:R\times R\rightarrow H[\lambda]$ be the $2$-cocycle induced by $s$. And we know that $\chi$ does not depend on the choice of $s$. Thus $\Phi(d_R,d_H)([\chi])$ does not depend on the choice of $s$. Then we can define a Wells map related to $\mathcal{E}$ as follows.

\begin{defi}
Let $\mathcal{E}:0 \to H \stackrel{\alpha}\to E \stackrel{p}\to R \to 0$ be an abelian extension of Leibniz conformal algebra $R$ by a $\cb[\partial]$-module $H$. The map $\mathcal{W}:{\rm g}(R,H)\rightarrow H^2(R,H)$,
$$\mathcal{W}(d_R,d_H)=\Phi(d_R,d_H)([\chi]),$$
for any $(d_R,d_H)\in{\rm g}(R,H)$, is called the {\rm Wells map} related to $\mathcal{E}$.
\end{defi}
More generally, we have the following theorem.

\begin{thm}\label{thm6.6}
Let $\mathcal{E}:0 \to H \stackrel{\alpha}\to E \stackrel{p}\to R \to 0$ be an abelian extension of Leibniz confoamal algebra $R$ by a $\cb[\partial]$-module $H$. A pair $(d_R,d_H)\in{\rm Der}^L(R)\times{\rm Der}^L(H)$ is extensible if and only if $(d_R,d_H)\in{\rm g}(R,H)$ and $\mathcal{W}(d_R,d_H)=0$.
\end{thm}
\begin{proof}
Let $s$ be a section of $\mathcal{E}$. Identify $H$ with a $\cb[\partial]$-submodule of $E$. Then there is a $\cb[\partial]$-module isomorphism $\tau:R\oplus H\rightarrow E, (a,h)\mapsto s(a)+h,$ for any $a\in R,h\in H$.

On the one hand, if $(d_R,d_H)\in{\rm g}(R,H)$ and $\mathcal{W}(d_R,d_H)=0$, there exists a $\cb[\partial]$-module homomorphism $f:R\rightarrow H$ such that $\Phi(d_R,d_H)(\chi)=\delta(f)$. Then for any $s(a)+h\in E$, we can define a $\cb[\partial]$-module homomorphism $d_E:E\rightarrow E$ by
$$d_E(s(a)+h)=f(a)+s(d_R(a))+d_H(h).$$
Obviously, $d_E$ satisfies
\begin{eqnarray*}
&&\ \ \ \ \ \ \ \ \ \ \ \ \ \ \ \ \ \ \ \ \ \ \ \ \ \ \ \ \ \ d_E(\alpha(h))=d_E(h)=d_H(h)=\alpha(d_H(h)), \\
&&p(d_E(s(a)+h))=p(f(a)+s(d_R(a))+d_H(h))=p(s(d_R(a)))=d_R(a)=d_R(p(s(a)+h)),
\end{eqnarray*}
for any $a\in R,h\in H$. Moreover, for any $s(a_1)+h_1,s(a_2)+h_2\in E$, we have
\begin{eqnarray*}
&&d_E([(s(a_1)+h_1)_\lambda(s(a_2)+h_2)])\\
&=&d_E([s(a_1)_\lambda s(a_2)]+[s(a_1)_\lambda h_2]+[{h_1}_{\lambda}s(a_2)])\\
&=&d_E\big(\chi_\lambda(a_1,a_2)+s([{a_1}_\lambda a_2])+[s(a_1)_\lambda h_2]+[{h_1}_{\lambda}s(a_2)]\big)\\
&=&d_H(\chi_\lambda(a_1,a_2))+f([{a_1}_\lambda a_2])+s(d_R([{a_1}_\lambda a_2]))+d_H([s(a_1)_\lambda h_2])+d_H([{h_1}_{\lambda}s(a_2)])\\
&=&d_H(\chi_\lambda(a_1,a_2))+f([{a_1}_\lambda a_2])+s([{d_R(a_1)}_\lambda a_2]+[{a_1}_\lambda d_R(a_2)])\\
&&+d_H([s(a_1)_\lambda h_2])+d_H([{h_1}_{\lambda}s(a_2)])\\
&=&d_H(\chi_\lambda(a_1,a_2))+f([{a_1}_\lambda a_2])+[{s(d_R(a_1))}_\lambda s(a_2)]-\chi_\lambda(d_R(a_1),a_2)\\
&&+[{s(a_1)}_\lambda s(d_R(a_2))]-\chi_\lambda(a_1,d_R(a_2))+d_H({a_1}\triangleright_\lambda h_2)+d_H({h_1}\triangleleft_{\lambda}{a_2}),
\end{eqnarray*}
and
\begin{eqnarray*}
&&[d_E(s(a_1)+h_1)_\lambda(s(a_2)+h_2)]+[(s(a_1)+h_1)_\lambda d_E(s(a_2)+h_2)]\\
&=&[(f(a_1)+s(d_R(a_1))+d_H(h_1))_\lambda(s(a_2)+h_2)]+[(s(a_1)+h_1)_\lambda (f(a_2)+s(d_R(a_2))+d_H(h_2))]\\
&=&[s(d_R(a_1))_\lambda s(a_2)]+[s(d_R(a_1))_\lambda h_2]+[f(a_1)_{\lambda}s(a_2)]+[d_H(h_1)_{\lambda}s(a_2)]\\
&&+[s(a_1)_\lambda s(d_R(a_2))]+[s(a_1)_\lambda f(a_2)]+[s(a_1)_\lambda d_H(h_2)]+[{h_1}_\lambda s(d_R(a_2))]\\
&=&[s(d_R(a_1))_\lambda s(a_2)]+d_R(a_1)\triangleright_\lambda h_2+[f(a_1)_{\lambda}s(a_2)]+d_H(h_1)\triangleleft_{\lambda}s(a_2)\\
&&+[s(a_1)_\lambda s(d_R(a_2))]+[s(a_1)_\lambda f(a_2)]+{a_1}\triangleright_\lambda d_H(h_2)+{h_1}\triangleleft_\lambda s(d_R(a_2)).
\end{eqnarray*}
Since $\Phi(d_R,d_H)(\chi)=\delta(f)$, that is
\begin{eqnarray*}
d_H(\chi_\lambda(a_1,a_2))-\chi_\lambda(d_R(a_1),a_2)-\chi_\lambda(a_1,d_R(a_2))=-f([{a_1}_\lambda a_2])+[s(a_1)_\lambda f(a_2)]+[f(a_1)_\lambda s(a_2)],
\end{eqnarray*}
for any $a_1,a_2\in R$, and $(d_R,d_H)\in{\rm g}(R,H)$, then
\begin{eqnarray*}
d_H({a_1}\triangleright_\lambda h_2)=d_R(a_1)\triangleright_\lambda h_2+{a_1}\triangleright_\lambda d_H(h_2), \ \ \ \
d_H({h_1}\triangleleft_{\lambda}{a_2})=d_H(h_1)\triangleleft_{\lambda}{a_2}+{h_1}\triangleleft_{\lambda}d_R(a_2).
\end{eqnarray*}
Combining the above two equations, we have
$$d_E([(s(a_1)+h_1)_\lambda(s(a_2)+h_2)])=[d_E(s(a_1)+h_1)_\lambda(s(a_2)+h_2)]+[(s(a_1)+h_1)_\lambda d_E(s(a_2)+h_2)],$$
which implies $d_E\in{\rm Der}^L(E)$. Hence the pair $(d_R,d_H)\in{\rm g}(R,H)$ is extensible.

On the other hand, if there is a left derivation $d_E\in{\rm Der}^L(E)$ such that $d_E\circ \alpha=\alpha\circ d_H$ and $d_R\circ p=p\circ d_E$, then $d_H=d_E|_H.$ Since for any $a\in R$, $p\big(d_E(s(a))-s(d_R(a))\big)=0$, we can define a $\cb[\partial]$-module homomorphism $f:R\rightarrow H$ by
$$f(a)=d_E(s(a))-s(d_R(a)).$$
Then for any $a\in R, h\in H$,
\begin{eqnarray*}
&&d_H(a\triangleright_\lambda h)=d_E(a\triangleright_\lambda h)=d_E([s(a)_\lambda h])=[d_E(s(a))_\lambda h]+[s(a)_\lambda d_E(h)]\\
&&=[d_E(s(a))_\lambda h]+a\triangleright_\lambda d_E(h)=[s(d_R(a))_\lambda h]+a\triangleright_\lambda d_E(h)=d_R(a)\triangleright_\lambda h+a\triangleright_\lambda d_E(h).
\end{eqnarray*}
Similarly, $d_H({h}\triangleleft_{\lambda}{a})=d_H(h)\triangleleft_{\lambda}{a}+{h}\triangleleft_{\lambda}d_R(a)$. Which implies $(d_R,d_H)\in{\rm g}(R,H)$. Finally, we prove that $\Phi(d_R,d_H)(\chi)=\delta(f)$. For any $a,b\in R$, we have
\begin{eqnarray*}
&&\Phi(d_R,d_H)(\chi)_\lambda(a,b)\\
&=&d_H(\chi_\lambda(a,b))-\chi_\lambda(d_R(a),b)-\chi_\lambda(a,d_R(b))\\
&=&d_E([s(a)_\lambda s(b)])-d_E\circ s([a_\lambda b])-[s(d_R(a))_\lambda s(b)]+s([d_R(a)_\lambda b])\\
&&-[s(a)_\lambda s(d_R(b))]+s([a_\lambda d_R(b)])\\
&=&d_E([s(a)_\lambda s(b)])-d_E\circ s([a_\lambda b])-[s(d_R(a))_\lambda s(b)]+s\circ d_R([a_\lambda b])-[s(a)_\lambda s(d_R(b))]\\
&=&d_E([s(a)_\lambda s(b)])-f([a_\lambda b])-[s(d_R(a))_\lambda s(b)]-[s(a)_\lambda s(d_R(b))]\\
&=&[d_E(s(a))_\lambda s(b)]+[s(a)_\lambda d_E(s(b))]-f([a_\lambda b])-[s(d_R(a))_\lambda s(b)]-[s(a)_\lambda s(d_R(b))]\\
&=&-f([a_\lambda b])+[f(a)_\lambda s(b)]+[s(a)_\lambda f(b)]=-f([a_\lambda b])+a\triangleright_\lambda f(b)+f(a)\triangleleft_\lambda b \\
&=&\delta(f)_\lambda(a,b).
\end{eqnarray*}
Hence $\mathcal{W}(d_R,d_H)=0$.
\end{proof}

Let $\mathcal{E}:0 \to H \stackrel{\alpha}\to E \stackrel{p}\to R \to 0$ be an abelian extension of Leibniz conformal algebra $R$ by a $\cb[\partial]$-module $H$. Denote
$${\rm Der}^L_H(E)=\{d_E\in{\rm Der}^L(E)\mid d_E(H)\subseteq H\}.$$
Then ${\rm Der}^L_H(E)$ is a Lie subalgebra of ${\rm Der}^L(E)$. $d_E\in{\rm Der}_H(E)$ induces two maps ${d_E}|_H$ and ${\hat{d}}_E=p\circ d_E\circ s$, where $s$ is a section of $\mathcal{E}$.
Hence we obtain a map
$$\kappa:{\rm Der}^L_H(E)\rightarrow{\rm Der}^L(R)\times{\rm Der}^L(H),\ \ \ \ \ \ \ d_E\mapsto({\hat{d}}_E,{d_E}|_H).$$

\begin{lem}
With the above notations, ${\rm Im}(\kappa)\subseteq{\rm g}(R,H)$ and $\kappa$ induces a Lie algebra homomorphism $\bar{\kappa}:{\rm Der}_H(E)\rightarrow{\rm g}(R,H)$.
\end{lem}
Considering the image and kernel of $\bar{\kappa}$, we have the following.

\begin{thm}
Let $\mathcal{E}:0 \to H \stackrel{\alpha}\to E \stackrel{p}\to R \to 0$ be an abelian extension of Leibniz conformal algebra $R$ by a $\cb[\partial]$-module $H$. There is an exact sequence of vector spaces
$$0 \to Z^1(R,H) \stackrel{\iota}\to {\rm Der}^L_H(E) \stackrel{\bar{\kappa}}\to {\rm g}(R,H) \stackrel{\mathcal{W}}\to H^2(R,H),$$
where $\iota$ is the inclusion map.
\end{thm}
\begin{proof}
First, we prove that ${\rm Ker}(\bar{\kappa})=Z^1(R,H)$. Any $d_E\in{\rm Ker}(\bar{\kappa})$, implies ${\hat{d}}_E=0$ and ${d_E}|_H=0$. For any $a\in R$, we have $p\circ d_E\circ s(a)=0,$ that is $d_E\circ s(a)\in H$. And for $d_E\in{\rm Ker}(\bar{\kappa})$, $d_E\circ s\in Z^1(R,H)$. Thus there is a map
$$\varphi:{\rm Ker}(\bar{\kappa})\rightarrow Z^1(R,H),\ \ \ \ \ \ \ d_E\mapsto d_E\circ s.$$
Now we prove that $\varphi$ is a bijection. If $d_E\in{\rm Ker}(\bar{\kappa})$, and $\varphi(d_E)=0$, for $a\in R, h\in H$, we have
$$d_E(s(a)+h)={d_E}|_H(h)+\varphi(d_E)(a)=0.$$
This implies $\varphi$ is injective. Conversely, for any $f\in Z^1(R,H)$, we define a $\cb[\partial]$-module homomorphism
$$d_E:E\rightarrow E,\ \ \ \ \ \ \ d_E(s(a)+h)=f(a),$$
for any $a\in R,h\in H$. Obviously, ${d_E}|_H=0$, ${\hat{d}}_E(a)=p(f(a))=0$ and $\varphi(d_E)(a)=d_E(s(a))=f(a),$ for $a\in R$. Moreover, $d_E\in{\rm Der}^L(E)$ and so $d_E\in {\rm Ker}(\bar{\kappa})$. Thus $\varphi$ is surjective.

Next, we need to show that the sequence is exact at ${\rm g}(R,H)$. For any $(d_R,d_H)\in{\rm Im}(\bar{\kappa})$, which implies $(d_R,d_H)$ is extensible, we obtain that  $\mathcal{W}(d_R,d_H)=0$. Then $ {\rm Im}(\bar{\kappa})\subseteq{\rm Ker}(\mathcal{W})$. Conversely, for any $(d_R,d_H)\in{\rm Ker}(\mathcal{W})$, there exists a map $f\in C^1(R,H)$ such that $\Phi(d_R,d_H)(\chi)=\delta(f)$. Define $d_E:E\rightarrow E$ by
$$d_E(s(a)+h)=f(a)+s(d_R(b))+d_H(h),$$
for any $a\in R,h\in H$. Similar to the proof of Theorem \ref{thm6.6}, we get $d_E\in{\rm Der}(E)$ and $\bar{\kappa}(d_E)=(d_R,d_H)$. Then ${\rm Ker}(\mathcal{W})\subseteq {\rm Im}(\bar{\kappa})$.
\end{proof}

\bigskip
\noindent
{\bf Acknowledgements. } This work was financially supported by National
Natural Science Foundation of China (No.12201182, 11801141).

 \end{document}